\documentclass[leqno,11pt]{amsart}
\usepackage{amsmath,amscd,amsthm,amsxtra}
\usepackage{epsfig,graphics}
\usepackage{amssymb,latexsym}
\usepackage{mathrsfs}
\usepackage{eucal}
\usepackage{hyperref}
\usepackage{pictexwd,dcpic}

\newtheorem{thm}{\bf Theorem}[section]

\newtheorem{prop}[thm]{\bf Proposition}
\newtheorem{cor}[thm]{\bf Corollary}
\newtheorem{lem}[thm]{\bf Lemma}
\newtheorem{rem}[thm]{\bf Remark}
\newtheorem{ex}[thm]{\bf Example}

\newcommand{\bs}{\boldsymbol}
\newcommand{\A}{\mathbb{A}}

\newcommand{\B}{\mathbf{B}}
\newcommand{\cB}{\mathbb{B}}
\newcommand{\W}{\mathcal{W}}
\newcommand{\cP}{\mathscr{P}}

\newcommand{\pf}{\noindent{\bfseries Proof. }}
\newcommand{\ov}{\overline}

\newcommand{\ba}{\bs{\rm a}}
\newcommand{\bb}{\bs{\rm b}}
\newcommand{\bc}{\bs{\rm c}}
\newcommand{\bd}{\bs{\rm d}}

\newcommand{\M}{{\mathcal{M}}}

\newcommand{\T}{\mathcal{T}}
\newcommand{\gl}{\mathfrak{gl}}
\newcommand{\Z}{\mathbb{Z}}

\newcommand{\te}{\widetilde{e}}
\newcommand{\tf}{\widetilde{f}}
\newcommand{\g}{\mathfrak{g}}
\newcommand{\td}{\widetilde}

\newcommand{\mf}{\mathfrak}

\newcommand{\nw}{^{\nwarrow}}
\newcommand{\se}{^{\searrow}}
\newcommand{\wh}{\widehat}

\numberwithin{equation}{section}

\begin{document}
\title[ ]
{RSK correspondence and classically irreducible Kirillov-Reshetikhin crystals}
\author{JAE-HOON KWON}
\address{Department of Mathematics \\ University of Seoul   \\  Seoul 130-743, Korea }
\email{jhkwon@uos.ac.kr }

\thanks{This work was  supported by Basic Science Research Program through the National Research Foundation of Korea (NRF) 
funded by the Ministry of  Education, Science and Technology (No. 2011-0006735).}

\begin{abstract}
We give a new  combinatorial model of the Kirillov-Reshetikhin
crystals of type $A_n^{(1)}$ in terms of non-negative integral
matrices based on the classical RSK algorithm, which has a simple description of  the affine crystal
structure  without using the
promotion operator. We have a similar
description of the Kirillov-Reshetikhin crystals  associated to exceptional nodes in the Dynkin diagrams of classical affine or non-exceptional affine type, which are called classically irreducible together with those of type $A_{n}^{(1)}$. %A uniform description of the affine crystal structure on these, so called classically irreducible,  Kirillov-Reshetikhin crystals is also obtained using a Lusztig's involution on an underlying classical crystal.
\end{abstract}

\maketitle

%\setcounter{tocdepth}{1}
%\tableofcontents

\section{Introduction}

The Robinson-Schensted-Knuth (simply RSK) correspondence is a weight preserving
bijection from the set $\M_{m \times n}$ of $m\times n$
non-negative integral matrices to the set $\T_{m\times n}$ of pairs
of semistandard Young tableaux of the same shape with entries from
$m$ and $n$ letters, respectively \cite{K}.

The RSK map $\kappa$ has nice representation theoretic
interpretations from a viewpoint of the Kashiwara's crystal base
theory  \cite{Kas91}. In \cite{La}, Lascoux show
that $\M_{m\times n}$  has a $\gl_{m}\oplus
\gl_{n}$-crystal structure and
$\kappa$ is an isomorphism of crystals, where one can define a
$\gl_{m}\oplus \gl_{n}$-crystal structure on $\T_{m\times n}$ in an
obvious way following \cite{KN}. As an application, a non-symmetric Cauchy kernel expansion into a sum of product of Demazure characters is obtained. 

In \cite{K09}, we show that
  $\kappa$ can be extended to an isomorphism of $\gl_{m+n}$-crystals. Here
$\M_{m\times n}$ or $\T_{m\times n}$ can be regarded as a crystal
associated to a generalized Verma module over $\gl_{m+n}$. As an
application, a weight generating function of plane
partitions in a bounded region is given as a Demazure character of
$\gl_{m+n}$. (See also \cite{K10} for an application of RSK to the
crystal base of a modified quantized enveloping algebra of type
$A_{+\infty}$ and $A_\infty$.)

The purpose of this paper is to study  the RSK correspondence
further in this direction and discuss  its connection with affine
crystals. We observe that
$\M_{r\times (n-r)}$ has a  natural affine crystal structure of type
$A_{n-1}^{(1)}$ for $n\geq 2$ and $1\leq r\leq n-1$ by \cite{K09} and the symmetry
of the Dynkin diagram of $A_{n-1}^{(1)}$. For $s\geq 1$, we let $\M_{r\times (n-r)}^s$ be the
set of matrices in $\M_{r\times (n-r)}$  such that the length of a
maximal decreasing subsequence of its row or column word is no more
than $s$. Then as the main result in this paper, we show (Theorem \ref{main 1}) that  as an
affine crystal of type $A_{n-1}^{(1)}$,
\begin{equation}\label{isomorphism}
\M_{r\times (n-r)}^s\otimes T_{s\omega_r} \simeq \B^{r,s},
\end{equation}
where  $\B^{r,s}$ is a perfect crystal  \cite{KMN2} with highest
weight $s\omega_r$ or the
rectangular partition $(s^r)$  as a classical $\gl_{n}$-crystal, and $T_{s\omega_r}=\{t_{s\omega_r}\}$
is a crystal with ${\rm wt}(t_{s\omega_r})=s\omega_r$,
$\varepsilon_i(t_{s\omega_r})=\varphi_i(t_{s\omega_r})=-\infty$ for
all $i$.  

To prove \eqref{isomorphism},  two RSK maps $\kappa^{\nw}$ and
$\kappa^{\se}$ are considered, which map a matrix in $\M_{r\times
(n-r)}^s$  to a pair of semistandard Young tableaux of normal and
anti-normal shape, respectively. They turn out to be the
projections of $\M_{r\times (n-r)}^s$ to a classical crystal of type
$A_{n-1}$ corresponding to maximal parabolic subalgebras obtained from $A_{n-1}^{(1)}$ by removing the simple roots $\alpha_0$ and  $\alpha_r$ respectively. These two RSK maps play an important role in proving the regularity of $\M_{r\times (n-r)}^s\otimes T_{s\omega_r}$ and constructing  the
isomorphism in \eqref{isomorphism}. Note that $\M_{r\times (n-r)}$ can be regarded as  a limit of
the crystals $\B^{r,s}\otimes T_{-s\omega_r}$ as $s$ goes to
infinity.\vskip 2mm

\begin{figure}
\includegraphics[width=9cm, height=14cm]{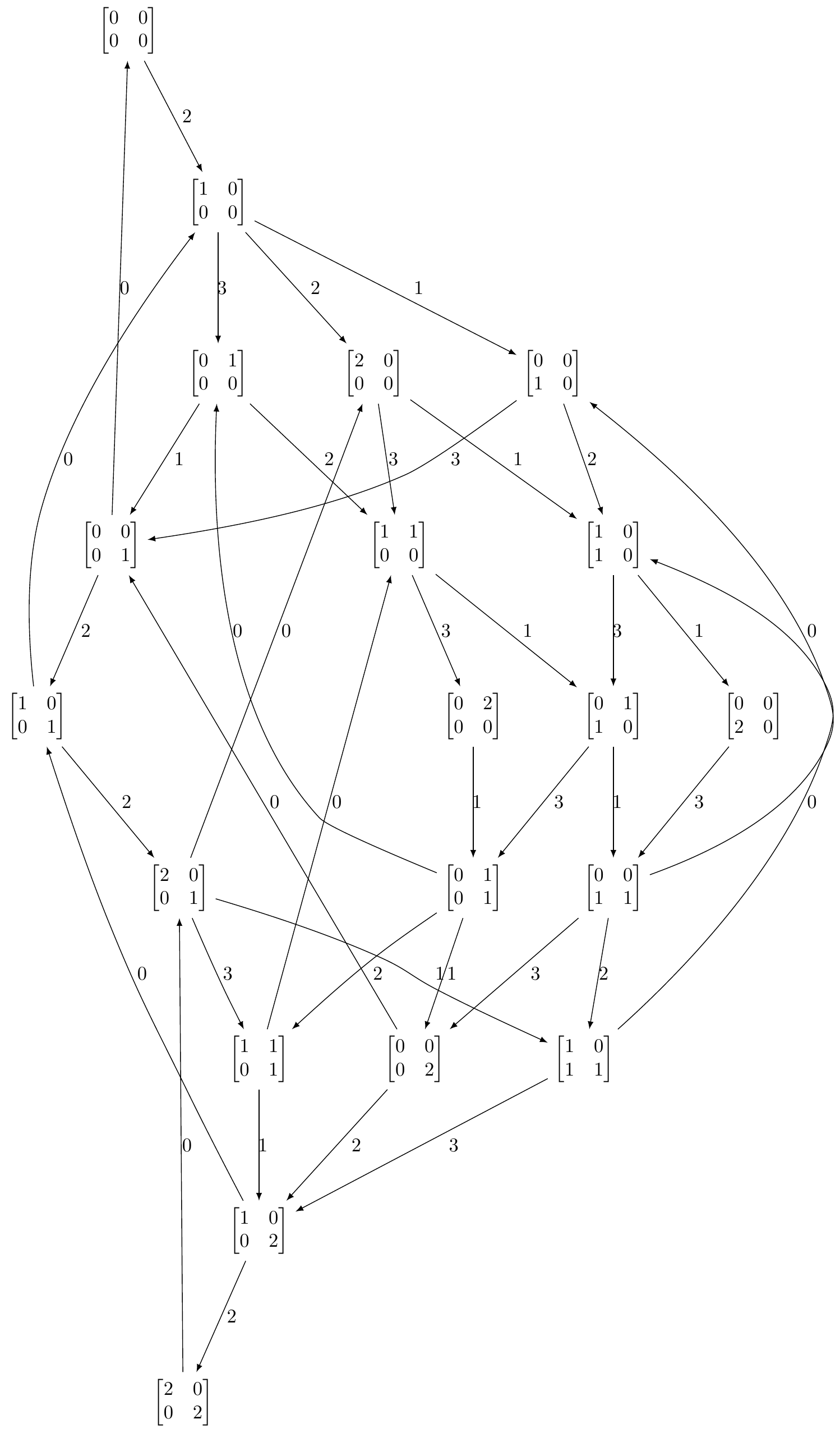}
\caption{The KR crystal $\B^{2,2}$ of type $A_3^{(1)}$ where the vertices are given in terms of non-negative integral $2\times 2$ matrices with the length of column or row words no more than 2. This graph was implemented by SAGE.} \label{Graph A}
\end{figure}

Let $\g$ be an affine  Kac-Moody algebra and let $U'_q(\g)$ be the
quantized enveloping algebra associated to the derived subalgebra
$\g'=[\g,\g]$. The finite dimensional irreducible $U_q'(\g)$-modules do not have crystal bases in general. But it was
conjectured by Hatayama et al. \cite{HKOTY,HKOTT} that a certain
family of finite dimensional irreducible $U_q'(\g)$-modules
$W^{r,s}$ called {\it Kirillov-Reshetikhin modules} (simply {\it KR modules}) \cite{KR}
have crystal bases $\B^{r,s}$, where  $r$ denotes a simple root
index of $\g$ except $0$ and $s$ is an arbitrary positive integer.
The conjectured crystals $\B^{r,s}$ are now called {\it KR crystals}.

For type $A_{n-1}^{(1)}$, the KR crystals $\B^{r,s}$ are the perfect
crystals in \eqref{isomorphism}. In this case, a
combinatorial description of $\B^{r,s}$ was given by Shimozono
\cite{S} using semistandard Young tableaux of  a rectangular shape
and the Sch\"{u}tzenberger's promotion operator \cite{Sc}. But, the main
advantage of our model using $r\times (n-r)$ integral matrices is that the description of its crystal
structure is remarkably simple, where the crystal operators or Kashiwara operators corresponding to $\alpha_0$ and $\alpha_r$ are given by adding  $\pm 1$ at the  entries at southeast and northwest corners of a matrix, respectively (see Figure ~\ref{Graph A}).

Recently, the existence of KR crystals $\B^{r,s}$ for the other classical affine
or non-exceptional affine  type was proved by Okado and Schilling \cite{OS},
and its combinatorial construction was given in \cite{FOS,OS}, where
the Kashiwara-Nakashima tableaux for the classical Lie algebras
\cite{KN} were used to describe the classical crystal structure on
$\B^{r,s}$.

We use \eqref{isomorphism} to obtain a new description of
the KR crystals associated to so-called {\it exceptional nodes} in the Dynkin diagrams of classical affine type (see Table 1 \cite{FOS}).  These crystals together with
$\B^{r,s}$ of type $A_{n-1}^{(1)}$ are called {\it classically
irreducible}  \cite{ST} since they are
connected as a classical crystal, and they are also perfect crystals \cite{KMN2}.

We use the Kashiwara's  method of folding crystals \cite{Kas96} to construct $\B^{n,s}$ of type $D_{n+1}^{(2)}$ and $C_n^{(1)}$ in terms of symmetric non-negative integral matrices (Theorem \ref{main 2}), and we describe $\B^{n-1,s}$ and $\B^{n,s}$ of type $D_n^{(1)}$ in terms of semi standard Young tableaux of type $A_{n-1}$ (Theorem \ref{main 3}). In both cases, the affine crystal structures are given explicitly as in $A_{n-1}^{(1)}$.

Finally we give some remarks on a relation between the affine crystal structure on classically irreducible KR crystals and the Lusztig involutions on an underlying classical crystal of type $A$.
\vskip 3mm

{\bf Acknowledgement} Part of this work was done during the author was visiting University of California, Berkeley in 2010-2011. The author would like to thank A. Schilling
for valuable discussion and kind explanation of KR crystals together with her recent preprint \cite{ST}.

\section{Preliminary}
\subsection{Quantum groups and crystals}
Let us give a brief review on crystals (cf. \cite{HK,Kas94}).
Let $A=(a_{ij})_{i,j\in I}$ be a generalized Cartan matrix with an index set $I$. 
Consider a quintuple
$(A,P^\vee,P,\Pi^\vee,\Pi)$ called a Cartan datum, where
\begin{itemize}
\item[]  $P^\vee$ is a free $\Z$-module of finite rank, 

\item[] $P={\rm Hom}_\Z(P^\vee,\Z)$,

\item[] $\Pi^\vee=\{\,h_i\,|\,i\in
I\,\}\subset P^\vee$,

\item[]   $\Pi=\{\,\alpha_i\,|\,i\in I\,\}\subset P$ such that $\langle\alpha_j,h_i\rangle=a_{ij}$ for $i,j\in I$.
\end{itemize}
A {\it crystal} associated to $(A,P^\vee,P,\Pi^\vee,\Pi)$ is a set
$B$ together with the maps ${\rm wt} : B \rightarrow P$,
$\varepsilon_i, \varphi_i: B \rightarrow \mathbb{Z}\cup\{-\infty\}$ and
$\te_i, \tf_i: B \rightarrow B\cup\{{\bf 0}\}$ ($i\in I$) such that
for $b\in B$ and $i\in I$
\begin{itemize}
\item[(1)] $\varphi_i(b) =\langle {\rm wt}(b),h_i \rangle +
\varepsilon_i(b),$

\item[(2)]  $\varepsilon_i(\te_i b) = \varepsilon_i(b) - 1$, $\varphi_i(\te_i b) =
\varphi_i(b) + 1$, ${\rm wt}(\te_ib)={\rm wt}(b)+\alpha_i$ if $\te_i b \neq {\bf 0}$,

\item[(3)] $\varepsilon_i(\tf_i b) = \varepsilon_i(b) + 1$, $\varphi_i(\tf_i b) =
\varphi_i(b) - 1$, ${\rm wt}({\tf_i}b)={\rm wt}(b)-\alpha_i$ if $\tf_i b \neq {\bf 0}$,

\item[(4)] $\tf_i b = b'$ if and only if $b = \te_i b'$ for $b, b' \in B$,

\item[(5)] $\te_ib=\tf_ib={\bf 0}$ if $\varphi_i(b)=-\infty$,
\end{itemize}
where ${\bf 0}$ is a formal symbol. Here we assume that
$-\infty+n=-\infty$ for all $n\in\Z$. Note that $B$ is equipped with
an $I$-colored oriented graph structure, where
$b\stackrel{i}{\rightarrow}b'$ if and only if $b'=\tf_{i}b$ for
$b,b'\in B$ and $i\in I$. We call $B$ { connected} if it is
connected as a graph, and  call  $B$ { normal} if
$\varepsilon_i(b)={\rm max}\{\,k\,|\,\te_i^kb\neq {\bf 0}\,\}$ and
$\varphi_i(b)={\rm max}\{\,k\,|\,\tf_i^kb\neq {\bf 0}\,\}$ for $b\in
B$ and $i\in I$. The { dual crystal $B^\vee$ of $B$} is defined
to be the set $\{\,b^\vee\,|\,b\in B\,\}$ with ${\rm
wt}(b^\vee)=-{\rm wt}(b)$, $\varepsilon_i(b^\vee)=\varphi_i(b)$,
$\varphi_i(b^\vee)=\varepsilon_i(b)$, $\te_i(b^\vee)=\left(\tf_i b
\right)^\vee$ and $\tf_i(b^\vee)=\left(\te_i b \right)^\vee$ for
$b\in B$ and $i\in I$. We assume that ${\bf 0}^\vee={\bf 0}$.

Let $B_1$ and $B_2$ be crystals. A {morphism}
$\psi : B_1 \rightarrow B_2$ is a map from $B_1\cup\{{\bf 0}\}$ to
$B_2\cup\{{\bf 0}\}$ such that
\begin{itemize}
\item[(1)] $\psi(\bf{0})=\bf{0}$,

\item[(2)] ${\rm wt}(\psi(b))={\rm wt}(b)$,
$\varepsilon_i(\psi(b))=\varepsilon_i(b)$, and
$\varphi_i(\psi(b))=\varphi_i(b)$ if $\psi(b)\neq \bf{0}$,

\item[(3)] $\psi(\te_i b)=\te_i\psi(b)$ if $\psi(b)\neq \bf{0}$ and
$\psi(\te_i b)\neq \bf{0}$,

\item[(4)] $\psi(\tf_i
b)=\tf_i\psi(b)$ if $\psi(b)\neq
\bf{0}$ and $\psi(\tf_i b)\neq \bf{0}$,
\end{itemize}
for $b\in B_1$ and $i\in I$.
We call $\psi$ an { embedding} and $B_1$ a { subcrystal of}
$B_2$ when $\psi$ is injective, and call $\psi$ {strict} if
$\psi : B_1\cup\{{\bf 0}\} \rightarrow B_2\cup\{{\bf 0}\}$ commutes
with $\te_i$ and $\tf_i$ for all $i\in I$, where we assume that $\te_i{\bf
0}=\tf_i{\bf 0}={\bf 0}$. When $\psi$ is a bijection, it is called an isomorphism.

For $b_i\in B_i$ ($i=1,2$), we say that { $b_1$ is equivalent to $b_2$} if
there exists an isomorphism of crystals $C(b_1)\rightarrow C(b_2)$
sending $b_1$ to $b_2$, where $C(b_i)$ is the connected component in $B_i$ including $b_i$ as an $I$-colored oriented graph.

A { tensor product $B_1\otimes B_2$} of crystals $B_1$ and $B_2$
is defined to be $B_1\times B_2$  as a set with elements  denoted by
$b_1\otimes b_2$, where  {\allowdisplaybreaks
\begin{equation*}
\begin{split}
{\rm wt}(b_1\otimes b_2)&={\rm wt}(b_1)+{\rm wt}(b_2), \\
\varepsilon_i(b_1\otimes b_2)&= {\rm
max}\{\varepsilon_i(b_1),\varepsilon_i(b_2)-\langle {\rm
wt}(b_1),h_i\rangle\}, \\
\varphi_i(b_1\otimes b_2)&= {\rm max}\{\varphi_i(b_1)+\langle {\rm
wt}(b_2),h_i\rangle,\varphi_i(b_2)\},\\
{\te}_i(b_1\otimes b_2)&=
\begin{cases}
{\te}_i b_1 \otimes b_2, & \text{if $\varphi_i(b_1)\geq \varepsilon_i(b_2)$}, \\
b_1\otimes {\te}_i b_2, & \text{if
$\varphi_i(b_1)<\varepsilon_i(b_2)$},
\end{cases}\\
{\tf}_i(b_1\otimes b_2)&=
\begin{cases}
{\tf}_i b_1 \otimes b_2, & \text{if  $\varphi_i(b_1)>\varepsilon_i(b_2)$}, \\
b_1\otimes {\tf}_i b_2, & \text{if $\varphi_i(b_1)\leq
\varepsilon_i(b_2)$},
\end{cases}
\end{split}
\end{equation*}
\noindent for $i\in I$. Here we assume that ${\bf 0}\otimes
b_2=b_1\otimes {\bf 0}={\bf 0}$.} Then $B_1\otimes B_2$ is a
crystal.

Let $\mathfrak{g}$ be a symmetrizable
Kac-Moody algebra associated  to
$A$. Let $P^\vee$ be  the dual weight lattice, $P={\rm
Hom}_\Z(P^\vee,\Z)$ the weight lattice, $\Pi^\vee=\{\,h_i\,|\,i\in
I\,\}$ the set of simple coroots, and
$\Pi=\{\,\alpha_i\,|\,i\in I\,\}$ the set of simple roots
of $\g$.
%Let $W$ be the Weyl group of $\g$, that is, the subgroup of $GL(P)$ generated
%by $r_i$ ($i\in I$), where $r_i$ is the simple reflection given by
%$r_i(\lambda)=\lambda-\langle \lambda,h_i \rangle\alpha_i$ for
%$\lambda\in P$.

Let $U_q(\g)$ be the quantized enveloping algebra of $\g$ over
$\mathbb{Q}(q)$  generated by $e_i$, $f_i$ and $q^h$ for $i\in I$
and $h\in P^\vee$. For a dominant integral weight $\Lambda$, let
$\B(\pm\Lambda)$ be the crystal base of an irreducible highest (resp.
lowest) weight $U_q(\g)$-module with highest (resp. lowest) weight
$\pm\Lambda$. Then $\B(\pm\Lambda)$ is a crystal associated to $(A,P^\vee,P,\Pi^\vee,\Pi)$.

We say that a crystal $B$ is {regular} if  $B$ is isomorphic to
the crystal base of an integrable $U_q(\g_J)$-module for any $J\subset
I$ with $|J|\leq 2$, where $\g_J$ is the Kac-Moody algebra associated with $A_J=(a_{ij})_{i,j\in J}$. Note that a regular crystal is normal.

% A regular crystal $B$ admits an action of the Weyl group on $B$ as follows; for
%$i\in I$ and $b\in B$
%\begin{equation}\label{Weyl}
%S_{r_i} b=
%\begin{cases}
%\tf_i^{\,\langle {\rm wt} (b),h_i \rangle}b, & \text{if $\langle {\rm
%wt} (b),h_i \rangle\geq 0$}, \\
%\te_i^{\,-\langle {\rm wt} (b),h_i \rangle}b, & \text{if $\langle {\rm
%wt} (b),h_i \rangle\leq 0$}, \\
%\end{cases}
%\end{equation}
%and $\eta=S_{r_{i_1}}\cdots S_{r_{i_t}}$ for  $w\in W$ with a reduced expression $w=r_{i_1}\cdots r_{i_t}$.

For $\Lambda\in P$, we denote by
$T_\Lambda=\{\,t_\Lambda\,\}$  a crystal with ${\rm
wt}(t_{\Lambda})=\Lambda$ and
$\varepsilon_i(t_{\Lambda})=\varphi_i(t_{\Lambda})=-\infty$ for $i\in I$.

\subsection{Quantum affine algebra}
Assume that $A$ is a generalized Cartan matrix
of affine type with an index set $I=\{\,0,1,\ldots,n\,\}$ following
\cite[\textsection 4.8]{K}, and  $\g$ is the associated affine Kac-Moody algebra  with the Cartan subalgebra ${\mf h}$.  Let $P^\vee=\bigoplus_{i\in I}\Z
h_i \oplus \Z d\subset {\mf h}$ be the dual weight lattice of $\g$, where $d$ is given by $\langle
\alpha_j,d \rangle=\delta_{0j}$ for $j\in I$. Let $\delta=\sum_{i\in I}a_i\alpha_i\in {\mf h}^*$ be the positive
imaginary null root of $\g$ and let $\Lambda_i\in {\mf h}^*$ ($i\in I$) be the $i^{\rm th}$ fundamental weight such that $\langle
\Lambda_i, h_j \rangle=\delta_{ij}$ for $j\in I$ and $\langle
\Lambda_i,d \rangle =0$. Then the weight lattice of $\g$ is $P=\bigoplus_{i\in I}\Z\Lambda_i\oplus\Z\frac{1}{a_0}\delta$.

Let $P_{\rm cl}=P/(\mathbb{Q}\delta\cap P)=\bigoplus_{i\in
I}\Z\Lambda_i$  and $(P_{\rm cl})^\vee=\bigoplus_{i\in I}\Z h_i$,
where  we still denote the image of $\Lambda_i$ in $P_{\rm cl}$  by
$\Lambda_i$.  Then we define $U'_q(\g)$ to be the subalgebra of
$U_q(\g)$ generated by $e_i$, $f_i$ and $q^h$ for $i\in I$ and $h\in
(P_{\rm cl})^\vee$. We regard $P_{\rm cl}$ as the weight lattice of
$U'_q(\g)$. For a proper subset $J \subset I$, let
$\Pi^\vee_J=\{\,h_i\,|\,i\in J\,\}$ and $\Pi_J=\{\,\alpha_i\,|\,i\in
J\,\}$, and let $U_q(\g_J)$ be the subalgebra of $U'_q(\g)$
generated by $e_i$, $f_i$ and $q^h$ for $i\in J$ and $h\in (P_{\rm
cl})^\vee$.

From now on, we mean by a $U'_q(\g)$-crystal (resp.
$U_q(\g_J)$-crystal)  a crystal associated to $(A,(P_{\rm
cl})^\vee,P_{\rm cl},\Pi^\vee,\Pi)$ (resp. $(A_J,(P_{\rm
cl})^\vee,P_{\rm cl},\Pi^\vee_J,\Pi_J)$). For simplicity, we will
often write the type of the generalized Cartan matrix $A$ (or $A_J$)
instead of $\g$ (or $\g_J$).

The following lemma plays an important role in this paper to have a combinatorial realization of KR crystals.
\begin{lem}[Lemma 2.6 in \cite{ST}]\label{ST lemma}
Let ${\mf g}$ be of classical affine or non-exceptional affine type. Fix $r\in I \setminus \{0\}$ and $s \geq 1$. Then any regular  $U'_q(\mf{g})$-crystal that is isomorphic to the KR crystal $\B^{r,s}$ as a $U_q(\mf{g}_{I \setminus \{0\}})$-crystal is also isomorphic to $\B^{r,s}$ as a $U'_q(\mf{g})$-crystal.
\end{lem}

\subsection{RSK algorithm}
Let us recall some necessary background on semistandard tableaux
following \cite{Fu,Mac95}. Let $\mathscr{P}$ be the set of
partitions. We identify a partition $\lambda=(\lambda_i)_{i\geq 1}$
with a {Young diagram} or a subset $\{\,(i,j)\,|\,1\leq j\leq
\lambda_i\,\}$ of $\mathbb{N}\times\mathbb{N}$. The {length of $\lambda$} is the number of
non-zero parts of $\lambda$ and  denoted by $\ell(\lambda)$. For $\mu\in\cP$ with $\mu\subset
\lambda$, $\lambda/\mu$ denotes the {skew Young diagram}, which
is  $\lambda\setminus\mu$ as a subset in
$\mathbb{N}\times\mathbb{N}$. We denote by $\lambda'=(\lambda'_i)_{i\geq 1}$
the { conjugate of $\lambda$}, and denote by $\lambda^\pi$ the
skew Young diagram obtained by $180^{\circ}$-rotation of $\lambda$.
For example,  \vskip 2mm
$$(5,3,2)\ \ = \ \
{\def\lr#1{\multicolumn{1}{|@{\hspace{.6ex}}c@{\hspace{.6ex}}|}{\raisebox{-.3ex}{$#1$}}}\raisebox{-.6ex}
{$\begin{array}{ccccc}
\cline{1-5}
\lr{\ \ } & \lr{\ \ } & \lr{\ \ } & \lr{\ \ }& \lr{\ \ }\\
\cline{1-5}
 \lr{\ \ } & \lr{\ \ }& \lr{\ \ } & & \\
\cline{1-3}
\lr{\ \ } & \lr{\ \ } & & & \\
\cline{1-2}
\end{array}$}} 
\ \ \ \ \ \ 
(5,3,2)^\pi \ \ = \ \
{\def\lr#1{\multicolumn{1}{|@{\hspace{.6ex}}c@{\hspace{.6ex}}|}{\raisebox{-.3ex}{$#1$}}}\raisebox{-.6ex}
{$\begin{array}{ccccc}
\cline{4-5}
& & &\lr{\ \ } & \lr{\ \ } \\
\cline{3-5}
& & \lr{\ \ } & \lr{\ \ }& \lr{\ \ }\\
\cline{1-5}
\lr{\ \ } & \lr{\ \ } & \lr{\ \ } & \lr{\ \ }& \lr{\ \ }\\
\cline{1-5}
\end{array}$}} \ \ \ .
$$\vskip 2mm

Let $\A$ be a linearly ordered set. For a skew Young diagram
$\lambda/\mu$,  let $SST_\A(\lambda/\mu)$ be the set of all
semistandard tableaux of shape $\lambda/\mu$ with entries in $\A$.
 Let $\W_{\A}$ be the set of finite words in
$\A$. We associate to each $T\in SST_\A(\lambda/\mu)$ a word
$w(T)\in\W_\A$, which is obtained by reading the entries of $T$ row
by row from top to bottom, and from right to left in each row.

Let ${\rm sh}(T)$ denote the shape of
a tableau $T$. If ${\rm sh}(T)=\nu$ (resp. $\nu^\pi$) for some $\nu\in\cP$,
then we say that $T$ is of {\it normal} (resp. {\it anti-normal})
shape. For $T\in SST_\A(\lambda/\mu)$, let $T^{\nw}$ (resp.
$T^{\se}$) be the unique semistandard tableau of normal (resp. anti normal)
shape such that $w(T^{\nw})$ (resp. $w(T^{\se})$) is Knuth equivalent to $w(T)$. Note that if ${\rm
sh}(T^{\nw})=\nu$, then ${\rm sh}(T^{\se})=\nu^\pi$.

Let $\mu\in\mathscr{P}$ and $a\in\A$ be given. For $S\in
SST_{\A}(\mu)$, let $a\rightarrow S$ be the tableau  obtained by
applying the Schensted's column insertion of $a$ into $S$. For
$w=w_1\ldots w_r\in\W_\A$, we define
\begin{equation*}
\begin{split}
{\bf P}(w)&=( w_r\rightarrow(\cdots(w_2\rightarrow w_1)\cdots)).
\end{split}
\end{equation*}

Let $\cB$ be another linearly ordered set. Let
\begin{equation}\label{MAB}
\M_{\A,\cB}=\left\{\,M=(m_{ab})_{a\in\A,b\in \cB}\,\,\Bigg\vert\,\, m_{ab}\in\Z_{\geq 0},\ \ \sum_{a,b}m_{ab}<\infty\,\right\}.
\end{equation}
Let $\Omega_{\A,\cB}$ be the set of biwords $(\ba,\bb)\in
\W_{\A}\times \W_{\cB}$ such that (1) $\ba=a_1\cdots a_r$ and
$\bb=b_1\cdots b_r$  for some $r\geq 0$, (2)   $(a_1,b_1)\leq \cdots
\leq (a_r,b_r)$, where for $(a,b)$ and $(c,d)\in \A\times \cB$,
$(a,b)< (c,d)$ if and only if $(b<d)$ or ($b=d$ and $a>c$). Then we have
a bijection from $\Omega_{\A,\cB}$ to $\M_{\A,\cB}$, where
$(\ba,\bb)$ is mapped to $M(\ba,\bb)=(m_{ab})$ with
$m_{ab}=\left|\{\,k\,|\,(a_k,b_k)=(a,b) \,\}\right|$. Note that the pair of
empty words $(\emptyset,\emptyset)$ corresponds to zero matrix. Let
$M\in \M_{\A,\cB}$ be given. Suppose that $M=M(\ba,\bb)$ and it transpose
$M^t=M(\bc,\bd)$ with $(\bc,\bd)\in\Omega_{\cB,\A}$. Let ${\bf P}(M)={\bf P}(\ba)$ and ${\bf Q}(M)={\bf
P}(\bc)$. Then we have a bijection called RSK correspondence:
\begin{equation*}
\kappa : \M_{\A,\cB}\ \ \longrightarrow\ \ \bigsqcup_{\lambda} SST_{\A}(\lambda)\times SST_{\cB}(\lambda),
\end{equation*}
where $M$ is mapped to $({\bf P}(M),{\bf Q}(M))$, and the union is over all $\lambda$ with $SST_{\A}(\lambda)\neq \emptyset$ and $SST_{\cB}(\lambda)\neq \emptyset$.

\section{KR crystals of type $A_{n-1}^{(1)}$}

\subsection{Affine algebra of type $A_{n-1}^{(1)}$}
Assume that $\g=A_{n-1}^{(1)}$ ($n\geq 2$) with
$I=\{\,0,1,\ldots,n-1\,\}$. We put $I_r=I\setminus\{r\}$ for $r\in
I$, and $I_{0,r}=I_0\cap I_r$ for $r\in I_0$. Note that
$\g_{I_0}\simeq\g_{I_r}=A_{n-1}$ and $\g_{I_{0,r}}=A_{r-1}\oplus
A_{n-r-1}$.

Let $\epsilon_k=\Lambda_k-\Lambda_{k-1}$ for $k\in I_0$ and
$\epsilon_n=\Lambda_0-\Lambda_{n-1}$.  Then
$\epsilon_1+\cdots+\epsilon_n=0$ and $\bigoplus_{i=1}^n\Z\epsilon_i$
forms a weight lattice of $\g_{I_0}$. Note that $\alpha_i=\epsilon_i -\epsilon_{i+1}$ for $i\in I_0$ and $\alpha_0=\epsilon_n-\epsilon_1$ in $P_{\rm cl}$. The fundamental weights for
$\g_{I_0}$ are
$\omega_i=\Lambda_i-\Lambda_0=\sum_{k=1}^i\epsilon_k$ for $i\in
I_0$.  

We regard $[n]=\{\,1<\cdots<n\,\}$ as a
$U_q(\g_{I_0})$-crystal $\B(\omega_1)$ with ${\rm
wt}(k)=\epsilon_k$, and
$[\ov{n}]=\{\,\ov{n}<\cdots< \ov{1}\,\}$ as its dual crystal with ${\rm wt}(\ov{k})=-\epsilon_k$. Then
$\W_{[n]}$ and $\W_{[\ov{n}]}$ are regular $U_q(\g_{I_0})$-crystals, where
we identify $w=w_1\ldots w_r$ with $w_1\otimes \cdots \otimes w_r$.

The fundamental weights for $\g_{I_r}$ are
$\omega'_i=\Lambda_i-\Lambda_{r}$ for $i\in I_r$. Note that 
$\omega_r=-\omega'_0$. In this case, we may identify a
$U_q(\g_{I_r})$-crystal $\B(\omega'_{r+1})$, the crystal base of the
natural representation of $U_q(\g_{I_r})$, with
$[n]_{+r}=\{\,r+1\prec \cdots\prec n\prec 1\prec \cdots\prec
r\,\}$.

\subsection{Affine crystal $\M_{r\times (n-r)}$}

For $1\leq r\leq n-1$, let
\begin{equation}
\M_{r\times (n-r)}=\M_{[\ov{r}], [n]\setminus [r]}
\end{equation}
(see \eqref{MAB}). First note that $\M_{r\times (n-r)}$ is a $U_q(A_{r-1})$-crystal with respect to
$\te_i$, $\tf_i$ ($1\leq i\leq r-1$),
where $\widetilde{x}_iM=M(\widetilde{x}_i\ba,\bb)$ for $x=e,f$ and $M\in
\M_{r\times (n-r)}$ with $M=M(\ba,\bb)$. Here, we assume that $\widetilde{x}_i M={\bf 0}$ if
$\widetilde{x}_i\ba={\bf 0}$. In a similar way, we may view
$\M_{r\times (n-r)}$ as a $U_q(A_{n-r-1})$-crystal with respect to
$\te_i$, $\tf_i$ ($r+1\leq i\leq n-1$) by considering the transpose
of $M\in \M_{r\times (n-r)}$ as an element in $\M_{[n]\setminus [r],[\ov{r}]}$. Since $\g_{I_{0,r}}\simeq
A_{r-1}\oplus A_{n-r-1}$,  $\M_{r\times (n-r)}$ is a regular
$U_q(\g_{I_{0,r}})$-crystal with 
\begin{equation}
{\rm wt}(M)=\sum_{i,j}m_{\ov{i}
j}(\epsilon_j-\epsilon_i),
\end{equation}
 for $M=(m_{\ov{i} j})\in
\M_{r\times (n-r)}$.

Now, let us define two more operators $\td{x}_0$ and $\td{x}_r$
($x=e,f$)  to make  $\M_{r\times (n-r)}$ a
$U'_q(A_{n-1}^{(1)})$-crystal. For $M=(m_{\ov{i} j})\in \M_{r\times
(n-r)}$, we define
\begin{equation}
\begin{split}
\te_r M&=
\begin{cases}
M-E_{\ov{r}\, r+1}, & \text{if $m_{\ov{r}\, r+1}\geq 1$}, \\
{\bf 0}, & \text{otherwise},
\end{cases}\ \ \ \ \
 \tf_r M= M+E_{\ov{r}\, r+1},\\
\tf_0 M&=
\begin{cases}
M-E_{\ov{1}\, n}, & \text{if $m_{\ov{1}\, n}\geq 1$}, \\
{\bf 0}, & \text{otherwise},
\end{cases}\ \ \ \ \ \ \ \ \ \ \
 \te_0 M= M+E_{\ov{1}\, n},
\end{split}
\end{equation}
where $E_{\ov{i}\, j}\in \M_{r\times (n-r)}$ denotes the elementary matrix with 1 at
the position  $(\ov{i},j)$ and $0$ elsewhere. Put
\begin{equation}
\begin{split}
\varepsilon_r(M)&=\max\left\{\,k\ \Big|\ \te_r^k M\neq {\bf 0}\,\right\},\ \
\varphi_r(M)=\varepsilon_r(M)+ \langle {\rm wt}(M), h_r \rangle,\\
\varphi_0(M)&=\max\left\{\,k\ \Big|\ \tf_0^k M\neq {\bf 0}\,\right\}, \ \ \varepsilon_0(M)=\varphi_0(M)- \langle {\rm wt}(M), h_0 \rangle.\\
\end{split}
\end{equation}
Then  we have
\begin{prop}\label{affine crystal M}
$\M_{r\times (n-r)}$ is a  $U'_q(A_{n-1}^{(1)})$-crystal with respect to ${\rm wt}$, $\varepsilon_i, \varphi_i$ and
$\te_i, \tf_i$ $(i\in I)$.
\end{prop}

\subsection{Young tableau descriptions of $\M_{r\times (n-r)}$ as a $U_q(A_{n-1})$-crystal}
Let us give another description of $\M_{r\times (n-r)}$  in terms of semistandard tableaux.
Consider
\begin{equation}
\begin{split}
\T_{r\times(n-r)}^{\se}&=\bigsqcup_{\ell(\lambda)\leq r,n-r}SST_{[\ov{r}]}(\lambda^\pi)\times SST_{[n]\setminus [r]}(\lambda^\pi).
\end{split}
\end{equation}
By \cite{KN}, $SST_{[\ov{r}]}(\lambda^\pi)\times
SST_{[n]\setminus[r]}(\lambda^\pi)$ is  a regular $U_q(\g_{I_{0,r}})$-crystal and
hence so is $\T_{r\times(n-r)}^{^{\searrow}}$. We will define operators
$\te_r, \tf_r$ on $\T_{r\times(n-r)}^{^{\searrow}}$ to make
$\T_{r\times(n-r)}^{^{\searrow}}$ a $U_q(\g_{I_{0}})$-crystal.

Let us first recall a combinatorial algorithm often called {\it a
signature rule}, which will be used throughout the paper. Suppose
that $\sigma=(\ldots,\sigma_{-2},\sigma_{-1},\sigma_0,\sigma_{1},\sigma_2,\ldots)$ be a sequence (not
necessarily finite) with $\sigma_{k}\in \{\,+\,,\,-\, , \ \cdot\ \}$ such that
$\sigma_k=\, +, $ or $\cdot\,$ for $k\gg 0$ and $\sigma_k=\, -$ or $\cdot$ for $k \ll 0$. In $\sigma$, we replace a pair
$(\sigma_{s},\sigma_{s'})=(+,-)$, where $s<s'$ and $\sigma_t=\,\cdot\,$
for $s<t<s'$, with $(\,\cdot\,,\,\cdot\,)$, and repeat this process
as far as possible until we get a sequence with no $-$ placed to the
right of $+$.  Such a reduced sequence will be denoted by $\td{\sigma}$.
When we have an infinite sequence $\sigma=(\sigma_1,\sigma_2,\ldots)$ (resp. $\sigma=(\ldots,\sigma_2,\sigma_1)$), we also understand $\tilde{\sigma}$ as a reduced sequence obtained by applying the signature rule to a  doubly infinite sequence $(\ldots,\,\cdot\,,\,\cdot\,,\,\cdot\,, \sigma_1,\sigma_2,\ldots)$ (resp. $(\ldots, \sigma_2,\sigma_1,\,\cdot\,,\,\cdot\,,\,\cdot\,)$).

Now, let $(S,T)\in \T_{r\times(n-r)}^{^{\searrow}}$ be given. For $k\geq 1$, let $s_k$
and $t_k$ be the entries in the top of the $k$-th columns of $S$ and
$T$ (enumerated from the right), respectively. We put
\begin{equation*}\label{signs}
\sigma_k=
\begin{cases}
+ \ , & \text{if the $k$-th column is empty}, \\
+ \ ,& \text{if $s_k> \ov{r}$ and $t_k>r+1$}, \\
- \ ,& \text{if $s_k=\ov{r}$ and $t_k=r+1$,}\\
\ \cdot \ \, ,& \text{otherwise}.
\end{cases}
\end{equation*}
Let $\td{\sigma}$ be the reduced sequence obtained from $\sigma=(\sigma_1,\sigma_2,\ldots)$ by the  signature rule.
Then we define $\te_r (S,T)$ to be the bitableaux obtained from
$(S,T)$  by removing $\boxed{\tiny{\ov{r}}}$ and  $\boxed{r+1}$ in
the columns of $S$ and $T$ corresponding to the right-most $-$ in
$\td{\sigma}$. If there is no such $-$ sign, then we define $\te_r
(S,T)={\bf 0}$. We define $\tf_r (S,T)$ to be the bitableaux
obtained from $(S,T)$ by adding $\boxed{\tiny{\ov{r}}}$ and  $\boxed{r+1}$
on top of the columns of $S$ and $T$ corresponding to the left-most
$+$ in $\td{\sigma}$. Note that $\tf_r^k(S,T)\neq {\bf 0}$ for all $k\geq 1$.

We put
\begin{equation}
\begin{split}
\varepsilon_r(S,T)&=\max\left\{\,k\ \Big|\ \te_r^k(S,T)\neq {\bf 0}\,\right\},\\
\varphi_r(S,T)&=\varepsilon_r(S,T)+ \langle {\rm wt}(S,T), h_r \rangle.\\
\end{split}
\end{equation}
Then $ \T_{r\times(n-r)}^{^{\searrow}}$ is a $U_q(\g_{I_0})$-crystal
with respect to ${\rm wt}$,  $\varepsilon_i, \varphi_i$ and $\te_i,
\tf_i$ ($i\in I_0$).

\begin{ex}\label{example of e,f on T_{m time n}}{\rm Suppose that $n=6$.
Consider
$$(S,T)\,=\left(\  \,{\def\lr#1{\multicolumn{1}{|@{\hspace{.6ex}}c@{\hspace{.6ex}}|}{\raisebox{-.3ex}{$#1$}}}\raisebox{-2.5ex}
{$\begin{array}[b]{cccc}
\cline{2-4}
 & \lr{\ov{3}} & \lr{\ov{2}}& \lr{\ov{2}}\\ 
\cline{1-4}
 \lr{\ov{3}} & \lr{\ov{2}} & \lr{\ov{1}}& \lr{\ov{1}}\\ 
\cline{1-4}
\end{array}$}} \ \  , \ \ 
\,{\def\lr#1{\multicolumn{1}{|@{\hspace{.6ex}}c@{\hspace{.6ex}}|}{\raisebox{-.3ex}{$#1$}}}\raisebox{-2.5ex}
{$\begin{array}[b]{cccc}
\cline{2-4}
 & \lr{{4}} & \lr{{4}}& \lr{4}\\ 
\cline{1-4}
 \lr{5} & \lr{5} & \lr{5}& \lr{6}\\ 
\cline{1-4}
\end{array}$}} \ \  \right).$$
Then
$$\te_3(S,T)\,=\left(\ \,{\def\lr#1{\multicolumn{1}{|@{\hspace{.6ex}}c@{\hspace{.6ex}}|}{\raisebox{-.3ex}{$#1$}}}\raisebox{-2.5ex}
{$\begin{array}[b]{cccc}
\cline{3-4}
 &   & \lr{\ov{2}}& \lr{\ov{2}}\\ 
\cline{1-4}
 \lr{\ov{3}} & \lr{\ov{2}} & \lr{\ov{1}}& \lr{\ov{1}}\\ 
\cline{1-4}
\end{array}$}} \ \ , \ \ 
\,{\def\lr#1{\multicolumn{1}{|@{\hspace{.6ex}}c@{\hspace{.6ex}}|}{\raisebox{-.3ex}{$#1$}}}\raisebox{-2.5ex}
{$\begin{array}[b]{cccc}
\cline{3-4}
 &   & \lr{{4}}& \lr{4}\\ 
\cline{1-4}
 \lr{5} & \lr{5} & \lr{5}& \lr{6}\\ 
\cline{1-4}
\end{array}$}} \ \  \right),$$
and
$$\tf_3(S,T)\,=\left(\  \,{\def\lr#1{\multicolumn{1}{|@{\hspace{.6ex}}c@{\hspace{.6ex}}|}{\raisebox{-.3ex}{$#1$}}}\raisebox{-2.5ex}
{$\begin{array}[b]{ccccc}
\cline{3-5}
& & \lr{\ov{3}} & \lr{\ov{2}}& \lr{\ov{2}}\\ 
\cline{1-5}
 \lr{\ov{3}} & \lr{\ov{3}} & \lr{\ov{2}} & \lr{\ov{1}}& \lr{\ov{1}}\\ 
\cline{1-5}
\end{array}$}} \ \  , \ \ 
\,{\def\lr#1{\multicolumn{1}{|@{\hspace{.6ex}}c@{\hspace{.6ex}}|}{\raisebox{-.3ex}{$#1$}}}\raisebox{-2.5ex}
{$\begin{array}[b]{ccccc}
\cline{3-5}
& & \lr{{4}} & \lr{{4}}& \lr{4}\\ 
\cline{1-5}
 \lr{4} & \lr{5} & \lr{5} & \lr{5}& \lr{6}\\ 
\cline{1-5}
\end{array}$}} \ \  \right).$$

}
\end{ex}

Define
\begin{equation}
\kappa^{^\searrow} : \M_{r\times (n-r)} \longrightarrow \T_{r\times(n-r)}^{^{\searrow}}
\end{equation}
by $\kappa^{^\searrow}(M) = ({\bf P}(M)^{^\searrow},
{\bf Q}(M)^{^\searrow})$.  By \cite[Theorem 3.6]{K09}, we have the
following.

\begin{prop}
$\kappa^{^\searrow}$ is an isomorphism of $U_q(\g_{I_0})$-crystals.
\end{prop}

\begin{ex}{\rm
Let $(S,T)$ be as in Example \ref{example of e,f on T_{m time n}}. Then $(S,T)=\kappa^{\se}(M)$, where $$M\,=\,
\begin{bmatrix}
1 & 0 & 1 \\
2 & 1 & 0 \\
0 & 2 & 0 \\
\end{bmatrix} .
$$
We have 
$$\te_3 M\,=\,
\begin{bmatrix}
0 & 0 & 1 \\
2 & 1 & 0 \\
0 & 2 & 0 \\
\end{bmatrix} 
$$
and $\kappa^{\se}(\te_3 M) =\te_3(S,T)$.
}
\end{ex}

Next, let us consider
\begin{equation}
\begin{split}
\T_{r\times(n-r)}^{^{\nwarrow}}&=\bigsqcup_{\ell(\lambda)\leq r, n-r}SST_{[\ov{r}]}(\lambda)\times SST_{[n]\setminus [r]}(\lambda).
\end{split}
\end{equation}
As in $\T_{r\times(n-r)}^{^{\searrow}}$,  $\T_{r\times(n-r)}^{^{\nwarrow}}$ is a regular
$U_q(\g_{I_{0,r}})$-crystal. Let us define operators $\te_0, \tf_0$ on
$\T_{r\times(n-r)}^{^{\nwarrow}}$ to make
$\T_{r\times(n-r)}^{^{\nwarrow}}$ a $U_q(\g_{I_{r}})$-crystal. Let
$(S,T)\in \T_{r\times(n-r)}^{^{\nwarrow}}$ be given.  For $k\geq 1$,
let $s_k$ and $t_k$ be the entries in the bottom of the $k$-th
columns of $S$ and $T$ (enumerated from the left), respectively. We
put
\begin{equation*}\label{signs}
\sigma_k=
\begin{cases}
- \ , & \text{if the $k$-th column is empty}, \\
- \ ,& \text{if $s_k<\ov{1}$ and $t_k<n$}, \\
+ \ ,& \text{if $s_k=\ov{1}$ and $t_k=n$,}\\
\ \cdot \ \, ,& \text{otherwise}.
\end{cases}
\end{equation*}
Let $\td{\sigma}$ be the reduced sequence obtained from $\sigma=(\ldots,\sigma_2,\sigma_1)$ by the signature rule.
We define $\te_0 (S,T)$ to be the bitableaux obtained from $(S,T)$ by adding $\boxed{\ov{1}}$ and
${\def\lr#1{\multicolumn{1}{|@{\hspace{.6ex}}c@{\hspace{.6ex}}|}{\raisebox{-.1ex}{$#1$}}}\raisebox{-.6ex}
{$\begin{array}[b]{c}
\cline{1-1}
\lr{n}\\
\cline{1-1}
\end{array}$}}$
below the bottom of the
columns of $S$ and $T$ corresponding to the right-most $-$ in $\td{\sigma}$. 
%If there is no such $+$ sign, then we define $\te_0 (S,T)={\bf 0}$. 
We define $\tf_0 (S,T)$ to be the bitableaux
obtained from $(S,T)$  by removing $\boxed{\ov{1}}$ and
${\def\lr#1{\multicolumn{1}{|@{\hspace{.6ex}}c@{\hspace{.6ex}}|}{\raisebox{-.1ex}{$#1$}}}\raisebox{-.6ex}
{$\begin{array}[b]{c}
\cline{1-1}
\lr{n}\\
\cline{1-1}
\end{array}$}}$
 in the columns of $S$ and $T$ corresponding to the left-most $+$ in $\td{\sigma}$.
If there is no such $+$ sign, then we define $\tf_0 (S,T)={\bf
0}$.  Note that $\te_0^k(S,T)\neq {\bf 0}$ for all $k\geq 1$.

We put
\begin{equation}
\begin{split}
\varphi_0(S,T)&=\max\left\{\,k\ \Big|\ \tf_0^k(S,T)\neq {\bf 0}\,\right\},\\
\varepsilon_0(S,T)&=\varphi_0(S,T)- \langle {\rm wt}(S,T), h_0 \rangle.\\
\end{split}
\end{equation}
Then $ \T_{r\times(n-r)}^{^{\nwarrow}}$ is a $U_q(\g_{I_r})$-crystal with respect to ${\rm wt}$, $\varepsilon_i, \varphi_i$ and
$\te_i, \tf_i$ ($i\in I_r$).
Define
\begin{equation}
\kappa^{^\nwarrow} : \M_{r\times (n-r)} \longrightarrow \T_{r\times(n-r)}^{^{\nwarrow}}
\end{equation}
by $\kappa^{^\nwarrow} (M) = ({\bf P}(M)^{\nw}, {\bf
Q}(M)^{\nw})=({\bf P}(M), {\bf
Q}(M))$. By the same argument as in \cite[Theorem 3.6]{K09},  we
have the following.

\begin{prop}
$\kappa^{^\nwarrow}$ is an isomorphism of $U_q(\g_{I_r})$-crystals.
\end{prop}

\subsection{Main Theorem}
For $M\in \M_{r\times (n-r)}$ with $M=M(\ba,\bb)$, let  $\ell(M)$ be the maximal length  of weakly decreasing subwords of $\ba$.
For $s\geq 1$, let
\begin{equation}
\M_{r\times (n-r)}^s=\{\,M\in \M_{r\times (n-r)}\,|\, \ell(M)\leq s \,\}.
\end{equation}
Note that $\ell(M)$ is the number of columns in ${\bf P}(M)$ or ${\bf Q}(M)$  (cf. \cite[\textsection 3.1]{Fu}).
We regard  $\M_{r\times (n-r)}^s$ as a subcrystal of $\M_{r\times (n-r)}$ and define a $U'_q(A_{n-1}^{(1)})$-crystal
\begin{equation}
\mathcal{B}^{r,s} =\M^s_{r\times (n-r)}\otimes T_{s\omega_r}.
\end{equation}

\begin{lem}\label{regularity of type A}
$\mathcal{B}^{r,s} $ is a regular $U'_q(\g)$-crystal that is isomorphic to $\B(s\omega_r)$ as a $U_q(\g_{I_0})$-crystal.
\end{lem}
\pf When restricted to $\M_{r\times (n-r)}^s$, we have the following bijections
\begin{equation}\label{restricted RSK}
\begin{split}
\kappa^{\se} : \M_{r\times (n-r)}^s &\longrightarrow \T_{r\times(n-r)}^{^{\searrow},s}, \\
\kappa^{\nw} : \M_{r\times (n-r)}^s &\longrightarrow \T_{r\times(n-r)}^{^{\nwarrow},s},
\end{split}
\end{equation}
where
\begin{equation*}
\begin{split}
 \T_{r\times(n-r)}^{^{\searrow},s}&=  \bigsqcup_{\substack{\ell(\lambda)\leq r, n-r \\ \lambda_1\leq s}}SST_{[\ov{r}]}(\lambda^\pi)\times SST_{[n]\setminus [r]}(\lambda^\pi), \\
 \T_{r\times(n-r)}^{^{\nwarrow},s}&=  \bigsqcup_{\substack{\ell(\lambda)\leq r, n-r \\ \lambda_1\leq s}}SST_{[\ov{r}]}(\lambda)\times SST_{[n]\setminus [r]}(\lambda).
\end{split}
\end{equation*}
Since $ \T_{r\times(n-r)}^{^{\searrow},s}$ (resp. $
\T_{r\times(n-r)}^{^{\nwarrow},s}$) can be viewed as a subcrystal of
$\T_{r\times(n-r)}^{^{\searrow}}$ (resp.
$\T_{r\times(n-r)}^{^{\nwarrow}}$), $\kappa^{\se}$ (resp. $\kappa^{\nw}$) is an isomorphism of $U_q(\g_{I_0})$ (resp.
$U_q(\g_{I_r})$)-crystals.

First we claim that $\T_{r\times(n-r)}^{^{\searrow},s}\otimes
T_{s\omega_r}$ is isomorphic to $\B(s\omega_r)$ as a
$U_q(\g_{I_0})$-crystal. Recall that $\B(s\omega_r)$ can be
identified with $SST_{[n]}((s^r))$ \cite{KN}.

Let $(S,T)\in \T_{r\times(n-r)}^{^{\searrow},s}$ be given where ${\rm sh}(S)={\rm sh}(T)=\lambda^\pi$ for some $\lambda\in\cP$ with $\lambda_1\leq s$.
Consider an isomorphism of  $U_q(\g_{\{1,\ldots,r-1\}})$-crystals, $$\varsigma : SST_{[\ov{r}]}(\lambda^\pi)\otimes T_{s\omega_r} \longrightarrow SST_{[r]}(\lambda^c),$$ where $\lambda^c=(s^r)\setminus \lambda^\pi =(s-\lambda_r,\ldots,s-\lambda_1)$ is a rectangular complement of $\lambda^\pi$ in $(s^r)$ (see \cite[Lemma 5.8]{K07} for an explicit description of $\varsigma$, which is given as $\sigma^s$). Let $S^c=\varsigma(S\otimes t_{s\omega_r})$ and let $U$ be the semistandard tableau in $SST_{[n]}((s^r))$ obtained by gluing $S^c$ and $T$. Therefore, the map sending $(S,T)\otimes t_{s\omega_r} $ to $U$ defines a weight preserving bijection (with the same notation)
\begin{equation}\label{covering 1}
\varsigma : \T_{r\times(n-r)}^{^{\searrow},s}\otimes T_{s\omega_r} \longrightarrow SST_{[n]}((s^r)).
\end{equation}
By definition, it is straightforward to check that $\varsigma$  commutes with $\te_r$ and $\tf_r$, which therefore implies that it is an isomorphism of $U_q(\g_{I_0})$-crystals.

Next consider $\T_{r\times(n-r)}^{^{\nwarrow},s}\otimes T_{s\omega_r}=\T_{r\times(n-r)}^{^{\nwarrow},s}\otimes T_{-s\omega'_0}$.  We claim that $\T_{r\times(n-r)}^{^{\nwarrow},s}\otimes T_{s\omega_r}$ is isomorphic to $\B(-s\omega'_0)$ as a $U_q(\g_{I_r})$-crystal. Since $\B(-s\omega'_0)=\B(s\omega'_{t})$ where $t\equiv 2r$ $\pmod{n}$,   $\B(-s\omega'_0)$ can be identified with $SST_{[n]_{+r}}((s^r))$.  

Let $(S,T)\in \T_{r\times(n-r)}^{^{\nwarrow},s}$ be given where ${\rm sh}(S)={\rm sh}(T)=\lambda$ for some $\lambda\in\cP$ with $\lambda_1\leq s$. By modifying the bijection in \cite[Lemma 5.8]{K07} (exchanging $k^\vee$ and $k$), we have an isomorphism of  $U_q(\g_{\{1,\ldots,r-1\}})$-crystals,
$$\ov{\varsigma} : SST_{[\ov{r}]}(\lambda)\otimes T_{s\omega_r} \longrightarrow SST_{[r]}((s^r)/ \lambda).$$
Let $\ov{S}^c=\ov{\varsigma}(S\otimes t_{s\omega_r})$ and let $U$ be the semistandard tableau in $SST_{[n]_{+r}}((s^r))$ obtained by gluing $\ov{S}^c$ and $T$. Then the map sending $(S,T)\otimes t_{s\omega_r} $ to $U$ defines a weight preserving bijection (with the same notation)
\begin{equation}\label{covering 2}
\ov{\varsigma} : \T_{r\times(n-r)}^{^{\nwarrow},s}\otimes T_{s\omega_r} \longrightarrow SST_{[n]_{+r}}((s^r)).
\end{equation}
As in (\ref{covering 1}), $\ov{\varsigma}$  commutes with $\te_0$ and $\tf_0$ and it is an isomorphism of $U_q(\g_{I_r})$-crystals.

Now, for a proper subset $J\subset I$ with $|J|\leq 2$, we have  $J\subset I_0$ or $J\subset I_r$ or $J\subset \{0,n\}$. By (\ref{covering 1}) and (\ref{covering 2}), $\mathcal{B}^{r,s} $ is a crystal base of an integrable $U_q(\g_J)$-module. Hence it is a regular $U'_q(A_{n-1}^{(1)})$-crystal.
\qed

\begin{ex}{\rm Assume that $n=6$, $r=3$.
Consider $$M\,=\,
\begin{bmatrix}
1 & 0 & 1 \\
2 & 1 & 0 \\
0 & 2 & 0 \\
\end{bmatrix}  \ \in \  \M^4_{3\times 3}.
$$
Then we have
$${\bf P}(M)^{^{\searrow}}\,=\,{\def\lr#1{\multicolumn{1}{|@{\hspace{.6ex}}c@{\hspace{.6ex}}|}{\raisebox{-.3ex}{$#1$}}}\raisebox{-2.5ex}
{$\begin{array}[b]{cccc}
\cline{2-4}
 & \lr{\ov{3}} & \lr{\ov{2}}& \lr{\ov{2}}\\ 
\cline{1-4}
 \lr{\ov{3}} & \lr{\ov{2}} & \lr{\ov{1}}& \lr{\ov{1}}\\ 
\cline{1-4}
\end{array}$}} \ \ \ , \ \ \ 
{\bf Q}(M)^{\searrow}\,=\,{\def\lr#1{\multicolumn{1}{|@{\hspace{.6ex}}c@{\hspace{.6ex}}|}{\raisebox{-.3ex}{$#1$}}}\raisebox{-2.5ex}
{$\begin{array}[b]{cccc}
\cline{2-4}
 & \lr{{4}} & \lr{{4}}& \lr{4}\\ 
\cline{1-4}
 \lr{5} & \lr{5} & \lr{5}& \lr{6}\\ 
\cline{1-4}
\end{array}$}} \ .$$
Note that as an element in a $U_q(A_2)$-crystal, ${\bf P}(M)^{\se}$ is equivalent to 
$$
{\def\lr#1{\multicolumn{1}{|@{\hspace{.6ex}}c@{\hspace{.6ex}}|}{\raisebox{-.3ex}{$#1$}}}\raisebox{-2.5ex}
{$\begin{array}[b]{cccc}
\cline{1-4}
\lr{1} & \lr{1} & \lr{3}& \lr{3}\\ 
\cline{1-4}
\lr{2} & & & \\ 
\cline{1-1}
\end{array}$}}\ .
$$
By gluing it with ${\bf Q}(M)^{\se}$, we have 
$$\ \ \ \ \ \ \ \ \ \ \ \ 
{\def\lr#1{\multicolumn{1}{|@{\hspace{.6ex}}c@{\hspace{.6ex}}|}{\raisebox{-.3ex}{$#1$}}}\raisebox{-2.5ex}
{$\begin{array}[b]{cccc}
\cline{1-4}
\lr{1} & \lr{1} & \lr{3}& \lr{3}\\ 
\cline{1-4}
\lr{2} & \lr{{4}} & \lr{{4}}& \lr{4}\\ 
\cline{1-4}
 \lr{5} & \lr{5} & \lr{5}& \lr{6}\\ 
\cline{1-4} 
\end{array}$}}\ \ \raisebox{2ex}{$ \in\B(4\omega_3)$}\ \ ,
$$
which is equivalent to $M\otimes t_{4\omega_3}\in \mathcal{B}^{3,4}$ as an element in a $U_q(\g_{I_0})$ ($=U_q(A_5)$)-crystal.
If we view $M\in \M^5_{4\times 3}$, then $M\otimes t_{5\omega_3}\in \mathcal{B}^{3,5}$ corresponds to 
$$\ \ \ \ \ \ \ \ \ \ \ \ 
{\def\lr#1{\multicolumn{1}{|@{\hspace{.6ex}}c@{\hspace{.6ex}}|}{\raisebox{-.3ex}{$#1$}}}\raisebox{-2.5ex}
{$\begin{array}[b]{ccccc}
\cline{1-5}
\lr{1} &\lr{1} & \lr{1} & \lr{3}& \lr{3}\\ 
\cline{1-5}
\lr{2} &\lr{2} & \lr{{4}} & \lr{{4}}& \lr{4}\\ 
\cline{1-5}
\lr{3} & \lr{5} & \lr{5} & \lr{5}& \lr{6}\\ 
\cline{1-5} 
\end{array}$}}\ \ \raisebox{2ex}{$ \in\B(5\omega_3)$}\ \ .
$$

On the other hand, we have
$${\bf P}(M)^{^{\nwarrow}}\,=\,{\def\lr#1{\multicolumn{1}{|@{\hspace{.6ex}}c@{\hspace{.6ex}}|}{\raisebox{-.3ex}{$#1$}}}\raisebox{-2.5ex}
{$\begin{array}[b]{cccc}
\cline{1-4}
\lr{\ov{3}} & \lr{\ov{3}} & \lr{\ov{2}}& \lr{\ov{2}}\\ 
\cline{1-4}
 \lr{\ov{2}} & \lr{\ov{1}}& \lr{\ov{1}} &\\ 
\cline{1-3}
\end{array}$}} \ \ \ , \ \ \ 
{\bf Q}(M)^{\nwarrow}\,=\,{\def\lr#1{\multicolumn{1}{|@{\hspace{.6ex}}c@{\hspace{.6ex}}|}{\raisebox{-.3ex}{$#1$}}}\raisebox{-2.5ex}
{$\begin{array}[b]{cccc}
\cline{1-4}
  \lr{{4}} & \lr{{4}}& \lr{4} & \lr{6}\\ 
\cline{1-4}
 \lr{5} & \lr{5} & \lr{5}& \\
\cline{1-3}
\end{array}$}} \ .$$
Note that as an element in a $U_q(A_2)$-crystal, ${\bf P}(M)^{\nw}$ is equivalent to 
$$
{\def\lr#1{\multicolumn{1}{|@{\hspace{.6ex}}c@{\hspace{.6ex}}|}{\raisebox{-.3ex}{$#1$}}}\raisebox{-2.5ex}
{$\begin{array}[b]{cccc}
\cline{4-4}
& & & \lr{1} \\ 
\cline{1-4}
\lr{1} & \lr{2} & \lr{3}& \lr{3}\\ 
\cline{1-4}
\end{array}$}}\ .
$$
By gluing it with ${\bf Q}(M)^{\nw}$, we have 
$$\ \ \ \ \ \ \ \ \ \ \ \ 
{\def\lr#1{\multicolumn{1}{|@{\hspace{.6ex}}c@{\hspace{.6ex}}|}{\raisebox{-.3ex}{$#1$}}}\raisebox{-2.5ex}
{$\begin{array}[b]{cccc}
\cline{1-4}
\lr{4} & \lr{4} & \lr{4}& \lr{6}\\ 
\cline{1-4}
\lr{5} & \lr{5} & \lr{5}& \lr{1}\\ 
\cline{1-4}
\lr{1} & \lr{2} & \lr{3}& \lr{3}\\ 
\cline{1-4} 
\end{array}$}}\ \ \raisebox{2ex}{$ \in\B(4\omega'_0)$}\ \ ,
$$
which is equivalent to $M\otimes t_{4\omega_3}\in \mathcal{B}^{3,4}$ as an element in a $U_q(\g_{I_3})$ ($=U_q(A_5)$)-crystal.
}
\end{ex}

\begin{thm}\label{main 1}
Let $\B^{r,s}$ be the KR crystal of type $A_{n-1}^{(1)}$ for $1\leq r\leq n-1$ and $s\geq 1$. Then as a $U'_q(A_{n-1}^{(1)})$-crystal, we have
$$\mathcal{B}^{r,s} \simeq \B^{r,s}.$$
\end{thm}
\pf  Note that $\B^{r,s}$ is isomorphic to $\B(s\omega_r)$ as a $U_q(\g_{I_0})$-crystal \cite{KMN2}. Then it follows from Lemmas \ref{ST lemma} and \ref{regularity of type A} that $\mathcal{B}^{r,s}\simeq \B^{r,s}$. \qed

\section{Classically irreducible KR crystals of type $D_{n+1}^{(2)}$ and $C_{n}^{(1)}$}
\subsection{Affine algebras of type $D_{n+1}^{(2)}$ and $C_{n}^{(1)}$}
Assume that $\g=A_{2n-1}^{(1)}$ ($n\geq 2$) with $I=\{\,0,1,\ldots,2n-1\,\}$ and the Cartan datum $(A,P^\vee,P,\Pi^\vee,\Pi)$, and  $\widehat{\g}=D_{n+1}^{(2)}$ or $C_n^{(1)}$ with $\widehat{I}=\{\,0,\ldots,n\,\}$ and the Cartan datum $(\widehat{A},\widehat{P}^\vee,\widehat{P},\widehat{\Pi}^\vee,\widehat{\Pi})$.

\begin{center}
\hskip -1cm \setlength{\unitlength}{0.16in} \medskip
\begin{picture}(17,6)
\put(2,2){\makebox(0,0)[c]{$A_{2n-1}^{(1)}\ :$}}

\put(14.25,-0.15){\makebox(0,0)[c]{$\bigcirc$}}
\put(14.25,4.15){\makebox(0,0)[c]{$\bigcirc$}}
\put(16.5,2){\makebox(0,0)[c]{$\bigcirc$}}

\put(14.6,0.2){\line(1,1){1.55}} \put(14.6,3.8){\line(1,-1){1.55}}

\put(8.35,4.15){\line(1,0){1.5}}
\put(12.28,4.15){\line(1,0){1.45}}
\put(11,4.15){\makebox(0,0)[c]{$\cdots$}}

\put(8.35,-0.15){\line(1,0){1.5}}
\put(12.28,-0.15){\line(1,0){1.45}}
\put(11,-0.15){\makebox(0,0)[c]{$\cdots$}}

\put(8,-0.15){\makebox(0,0)[c]{$\bigcirc$}}
\put(8,4.15){\makebox(0,0)[c]{$\bigcirc$}}
\put(5.7,2){\makebox(0,0)[c]{$\bigcirc$}}

\put(6,2.3){\line(1,1){1.62}} \put(6,1.7){\line(1,-1){1.62}}

\put(5.6,1){\makebox(0,0)[c]{\tiny ${\alpha}_0$}}
\put(8,-1){\makebox(0,0)[c]{\tiny ${\alpha}_{2n-1}$}}
\put(8,3.3){\makebox(0,0)[c]{\tiny ${\alpha}_1$}}

\put(14.2,3.3){\makebox(0,0)[c]{\tiny ${\alpha}_{n-1}$}}
\put(14.2,-1){\makebox(0,0)[c]{\tiny ${\alpha}_{n+1}$}}
\put(16.8,1){\makebox(0,0)[c]{\tiny ${\alpha}_{n}$}}
\end{picture}
\end{center}

\begin{center}
\hskip -1cm  \setlength{\unitlength}{0.16in}
\begin{picture}(17,4)
\put(2,2){\makebox(0,0)[c]{$D_{n+1}^{(2)}\ :$}}

\put(5.6,2){\makebox(0,0)[c]{$\bigcirc$}}
\put(8,2){\makebox(0,0)[c]{$\bigcirc$}}
\put(14.25,2){\makebox(0,0)[c]{$\bigcirc$}}
\put(16.5,2){\makebox(0,0)[c]{$\bigcirc$}}

\put(8.35,2){\line(1,0){1.5}}
\put(12.28,2){\line(1,0){1.45}}
\put(15.4,2){\makebox(0,0)[c]{$\Longrightarrow$}}

\put(6.8,1.97){\makebox(0,0)[c]{$\Longleftarrow$}}
\put(11,1.95){\makebox(0,0)[c]{$\cdots$}}

\put(5.6,1){\makebox(0,0)[c]{\tiny $\widehat{\alpha}_0$}}
\put(8,1){\makebox(0,0)[c]{\tiny $\widehat{\alpha}_1$}}

\put(14.4,1){\makebox(0,0)[c]{\tiny $\widehat{\alpha}_{n-1}$}}
\put(16.8,1){\makebox(0,0)[c]{\tiny $\widehat{\alpha}_{n}$}}
\end{picture}
\end{center}
\begin{center}
\hskip -1cm  \setlength{\unitlength}{0.16in}
\begin{picture}(17,4)
\put(2,2){\makebox(0,0)[c]{$C_{n}^{(1)}\ :$}}

\put(5.6,2){\makebox(0,0)[c]{$\bigcirc$}}
\put(8,2){\makebox(0,0)[c]{$\bigcirc$}}
%\put(10.4,2){\makebox(0,0)[c]{$\bigcirc$}}
%\put(14.85,2){\makebox(0,0)[c]{$\bigcirc$}}
\put(14.25,2){\makebox(0,0)[c]{$\bigcirc$}}
\put(16.5,2){\makebox(0,0)[c]{$\bigcirc$}}
\put(8.35,2){\line(1,0){1.5}} %\put(10.82,2){\line(1,0){0.8}}
%\put(13.2,2){\line(1,0){1.2}}
\put(12.28,2){\line(1,0){1.45}}
\put(15.4,2){\makebox(0,0)[c]{$\Longleftarrow$}}

\put(6.8,1.97){\makebox(0,0)[c]{$\Longrightarrow$}}
\put(11,1.95){\makebox(0,0)[c]{$\cdots$}}

\put(5.6,1){\makebox(0,0)[c]{\tiny $\widehat{\alpha}_0$}}
\put(8,1){\makebox(0,0)[c]{\tiny $\widehat{\alpha}_1$}}
%\put(10.4,1){\makebox(0,0)[c]{\tiny $\widehat{\alpha}_2$}}
%\put(15,1){\makebox(0,0)[c]{\tiny $\widehat{\alpha}_{n-2}$}}
\put(14.4,1){\makebox(0,0)[c]{\tiny $\widehat{\alpha}_{n-1}$}}
\put(16.8,1){\makebox(0,0)[c]{\tiny $\widehat{\alpha}_{n}$}}
\end{picture}
\end{center}

Throughout this section, we assume that $\epsilon\in \{1,2\}$ and $\wh{\g}=D_{n+1}^{(2)}$ (resp. $\wh{\g}=C_{n}^{(1)}$) when $\epsilon=1$ (resp. $\epsilon=2$).

Put $\wh{I}_{r}=\wh{I}\setminus \{r\}$ ($r=0,n$) and $\wh{I}_{0,n}=\wh{I}_{0}\cap \wh{I}_{n}$.  Note that $\wh{\g}_{\wh{I}_0}\simeq\wh{\g}_{\wh{I}_n}=B_{n}$ (resp. $C_n$) when $\epsilon=1$ (resp. $\epsilon=2$) and $\wh{\g}_{\wh{I}_{0,n}}\simeq A_{n-1}$.
 We may assume that
\begin{equation*}
\begin{split}
&\wh{P}^\vee=\Z h_0\oplus\cdots\oplus\Z h_n \oplus\Z d \subset P^\vee,\\
&\wh{P}=\{\,\lambda\,|\,  \tfrac{1}{\epsilon}\langle\lambda,h_i\rangle \in\Z \ \ (i=0,n),\ \  \langle\lambda,h_i\rangle=\langle\lambda,h_{2n-i}\rangle\ \ (i\in \wh{I}_{0,n})\,\} \subset P,\\
&\wh{\Pi}^\vee= \{\,  \widehat{h}_i=h_i \ (i\in \wh{I}) \,\}\subset \Pi^\vee, \\
&\wh{\Pi}=\{\,\widehat{\alpha}_i= \epsilon\alpha_i \ \ (i=0,n),\ \
\widehat{\alpha}_i=\alpha_{i}+\alpha_{2n-i} \ \ (i\in \wh{I}_{0,n})\,\}\subset \Pi.\\
\end{split}
\end{equation*}
%Since $\wh{\delta}=\wh{\alpha}_0+2\wh{\alpha}_{1}+\cdots+2\wh{\alpha}_{n-1}+\wh{\alpha}_n=\delta$ is the positive imaginary null root of $\wh{\g}$ with minimal coefficients,
The classical weight lattice of $\wh{\g}$ is  $\wh{P}_{\rm cl}=\bigoplus_{i\in \wh{I}}\Z\wh{\Lambda}_i$ and its dual classical weight lattice is $(\wh{P}_{\rm cl})^\vee=\bigoplus_{i\in \wh{I}}\Z h_i$, where $\wh{\Lambda}_i=\epsilon\Lambda_i$ for $i=0,n$ and $\wh{\Lambda}_i=\Lambda_i+\Lambda_{2n-i}$ for $i\in \wh{I}_{0,n}$.
Note that $\wh{\alpha}_i=\wh{\epsilon}_i-\wh{\epsilon}_{i+1}$ ($i\in I_{0,n}$), where $\wh{\epsilon}_i=\epsilon_i-\epsilon_{2n-i+1}$ for $i=1,\ldots, n$, $\wh{\alpha}_0=-\epsilon\,\wh{\epsilon}_1$ and $\wh{\alpha}_n=\epsilon\,\wh{\epsilon}_n$ in $\wh{P}_{\rm cl}$.

We denote the fundamental weights for $\wh{\g}_{I_0}$ by $\wh{\omega}_i=\omega_i+\omega_{2n-i}$ for $i\in \wh{I}_{0,n}$ and $\wh{\omega}_n=\epsilon \omega_n$, and  those for $\wh{\g}_{I_n}$ by $\wh{\omega}'_i=\omega'_i+\omega'_{2n-i}$ for $i\in \wh{I}_{0,n}$ and $\wh{\omega}'_0=\epsilon \omega'_0=-\wh{\omega}_n$.

\subsection{Crystals of symmetric matrices}
Put
\begin{equation}
\wh{\M}_n = \left\{\,M=(m_{\ov{i} j})\in \M_{n\times n}\ \Big|\,  \text{$m_{\ov{i} j}=m_{\ov{j} i}$ and $\epsilon \, |\, m_{\ov{i} i}$ for $i,j\in [n]$}\,\right\}.
\end{equation}
Define
\begin{equation}\label{folded operators}
\begin{split}
\wh{e}_i=
\begin{cases}
(\te_i)^\epsilon , & \text{for $i=0,n$}, \\
\te_i\te_{2n-i}, & \text{for $i\in \wh{I}_{0,n}$},
\end{cases}
\ \ \ \ \ 
\wh{f}_i =
\begin{cases}
\left(\tf_i\right)^\epsilon, & \text{for $i=0,n$}, \\
\tf_i\tf_{2n-i}, & \text{for $i\in \wh{I}_{0,n}$}.
\end{cases}
\end{split}
\end{equation}
Note that $\M_{n\times n}$ is a $U'_q(A_{2n-1}^{(1)})$-crystal with respect to  ${\rm wt}$, $\varepsilon_i, \varphi_i$ and $\te_i, \tf_i$ $(i\in I)$ by Proposition \ref{affine crystal M}. Then it is not difficult to see that $\wh{\M}_n \cup\{{\bf 0}\}$ is invariant under $\wh{e}_i$ and $\wh{f}_i$ for $i\in \wh{I}$ (cf. \cite[Proposition 5.14]{K09}).  For $M\in \wh{\M}_n$, define 
\begin{equation}
\begin{split}
\wh{\rm wt}(M)=&{\rm wt}(M), \\
\wh{\varepsilon}_i (M)=
\begin{cases}
\varepsilon_i(M), & \text{if $i=0,n$}, \\
\tfrac{1}{\epsilon}\varepsilon_0(M), & \text{if $i\in \wh{I}_{0,n}$}.
\end{cases}, & \ \ \ \ \
\wh{\varphi}_i (M)=
\begin{cases}
\varphi_i(M), & \text{if $i=0,n$}, \\
\tfrac{1}{\epsilon}\varphi_0(M), & \text{if $i\in \wh{I}_{0,n}$}.
\end{cases}
\end{split}
\end{equation}

Hence we have

\begin{prop}
$\wh{\M}_n$ is a $U'_q(\wh{\g})$-crystal with respect to  $\wh{\rm wt}$, $\wh{\varepsilon}_i, \wh{\varphi}_i$ and
$\wh{e}_i, \wh{f}_i$ $(i\in \wh{I})$.
\end{prop}

Next, consider
\begin{equation}
\begin{split}
\wh{\T}_{n}^{\se}=\bigsqcup_{\ell(\lambda)\leq n} SST_{[\ov{n}]}(\epsilon\lambda^\pi), \ \ \ \ \ \wh{\T}_{n}^{^{\nwarrow}}=\bigsqcup_{\ell(\lambda)\leq n} SST_{[\ov{n}]}(\epsilon\lambda),
\end{split}
\end{equation}
where $2\lambda=(2\lambda_i)_{i\geq 1}$ for $\lambda=(\lambda_i)_{i\geq 1}\in \cP$.
They are regular $U_q(\wh{\g}_{\wh{I}_{0,n}})$-crystals with respect to $\te_i$, $\tf_i$ ($i\in I_{0,n}$). Here 
\begin{equation}
{\rm wt}(T)=-\sum_{i\in [n]}m_{\ov{i}}\wh{\epsilon}_i,
\end{equation}
for $T\in \wh{\T}_{n}^{\se}$ or $\wh{\T}_{n}^{\nw}$, where $m_{\ov{i}}$ is the number of $\ov{i}$'s appearing in $T$. 

Let us define operators $\te_n$, $\tf_n$ on $\wh{\T}_{n}^{\se}$ corresponding to $\wh{\alpha}_n$ as follows:

\noindent \textsc{Case 1.} Suppose that $\epsilon=1$.  Let $T\in \wh{\T}_{n}^{\se}$ be given. For $k\geq 1$, let $t_k$  be the entry in the top of the $k$-th column of $T$ (enumerated from the right). Consider  $\sigma=(\sigma_1,\sigma_2,\ldots)$, where
\begin{equation*}
\sigma_k=
\begin{cases}
+ \ , & \text{if $t_k>\ov{n}$ or the $k$-th column is empty}, \\
- \ ,& \text{if   $t_k=\ov{n}$}.
\end{cases}
\end{equation*}
Then we define $\te_n T$ to be the tableau
obtained from $T$  by removing  $\boxed{\ov{n}}$ in the column  corresponding to the right-most $-$ in $\td{\sigma}$.  If there is no such $-$ sign, then we define $\te_n T={\bf 0}$.
We define $\tf_n T$ to be the tableau  obtained from $T$ by adding   $\boxed{\ov{n}}$ on top of the
column  corresponding to the left-most $+$ in $\td{\sigma}$. 
%If there is no such $+$ sign, then we define $\wh{f}_n T={\bf 0}$.
\vskip 2mm

\noindent\textsc{Case 2.} Suppose that $\epsilon=2$.  Let $T\in \wh{\T}_{n}^{\se}$ be given.
For each $k\geq 1$, let $(t_{2k},t_{2k-1})$  the pair of entries  in the top of the $2k$-th  and $(2k-1)$-st columns of $T$ (from the right), respectively. Note that  $t_{2k}$ and $t_{2k-1}$ are placed in the same row and $t_{2k}\leq t_{2k-1}$. Consider $\sigma=(\sigma_1,\sigma_2,\ldots)$, where
\begin{equation*}\label{signs}
\sigma_k=
\begin{cases}
+ \ , & \text{if  $t_{2k},t_{2k-1}> \ov{n}$ or the $(2k-1)$-st column is empty}, \\
- \ ,& \text{if $t_{2k}=t_{2k-1}=\ov{n}$}, \\
\ \cdot \ \, ,& \text{otherwise}.
\end{cases}
\end{equation*}
Then we define $\te_n T$ to be the tableau
obtained from $T$  by removing  a domino
${\def\lr#1{\multicolumn{1}{|@{\hspace{.6ex}}c@{\hspace{.6ex}}|}{\raisebox{-.3ex}{$#1$}}}\raisebox{-.6ex}
{$\begin{array}[b]{cc}
\cline{1-1}\cline{2-2}
\lr{\ov{n}}&\lr{\ov{n}}\\
\cline{1-1}\cline{2-2}
\end{array}$}}$
 in the pair of columns  corresponding to the right-most $-$ in $\td{\sigma}$.
If there is no such $-$ sign, then we define $\te_nT={\bf
0}$. We define $\tf_n T$ to be the tableau obtained from $T$ by adding  a domino
${\def\lr#1{\multicolumn{1}{|@{\hspace{.6ex}}c@{\hspace{.6ex}}|}{\raisebox{-.3ex}{$#1$}}}\raisebox{-.6ex}
{$\begin{array}[b]{cc}
\cline{1-1}\cline{2-2}
\lr{\ov{n}}&\lr{\ov{n}}\\
\cline{1-1}\cline{2-2}
\end{array}$}}$
 on top of the pair of
columns  corresponding to the left-most $+$ in $\td{\sigma}$. 
%If there is no such $+$ sign, then we define $\wh{f}_n T={\bf 0}$.

Hence $\wh{\T}_{n}^{\se}$ is a $U'_q(\wh{\g}_{\wh{I}_0})$-crystal with respect to ${\rm wt}$, $\varepsilon_i$, $\varphi_i$, $\te_i$, $\tf_i$ ($i\in \wh{I}_0$), where
\begin{equation}
\begin{split}
\varepsilon_n(T)&=\max\left\{\,k\ \Big|\ \te_n^kT\neq {\bf 0}\,\right\},\ \
\varphi_n(T)=\varepsilon_n(T)+ \langle {\rm wt}(T), \wh{h}_n \rangle.\\
\end{split}
\end{equation}

Define
\begin{equation}
\wh{\kappa}^{\se} : \wh{\M}_n \longrightarrow \wh{\T}_{n}^{\se},
\end{equation}
by $\wh{\kappa}^{\se}(M)={\bf P}(M)^{\se}$. By \cite[Propositions 3.5 and 6.5]{K11}, we have
\begin{prop}\label{symmetric kappa se}
$\wh{\kappa}^{\se}$ is an isomorphism of $U_q({\wh{\g}_{\wh{I}_0}})$-crystals.
\end{prop}

Let us define operators $\te_0$, $\tf_0$ on $\wh{\T}_{n}^{\nw}$ corresponding to $\wh{\alpha}_0$ as follows:

\noindent \textsc{Case 1.} Suppose that $\epsilon=1$.
Let  $T\in \wh{\T}_{n}^{\nw}$ be given.
For $k\geq 1$, let $t_k$  be the entry
in the bottom of the $k$-th column of $T$ (enumerated from the left). Consider  $\sigma=(\ldots,\sigma_2,\sigma_1)$, where
\begin{equation*}
\sigma_k=
\begin{cases}
- \ , & \text{if $t_k< \ov{1}$ or the $k$-th column is empty}, \\
+ \ ,& \text{if   $t_k=\ov{1}$}.
\end{cases}
\end{equation*}
Then we define $\te_0 T$ to be the tableau
obtained from $T$  by adding
${\def\lr#1{\multicolumn{1}{|@{\hspace{.6ex}}c@{\hspace{.6ex}}|}{\raisebox{-.1ex}{$#1$}}}\raisebox{-.6ex}
{$\begin{array}[b]{c}
\cline{1-1}
\lr{\ov{1}}\\
\cline{1-1}
\end{array}$}}$ to the bottom of the column  corresponding to the right-most $-$ in $\td{\sigma}$.  
We define $\tf_0 T$ to be the tableau  obtained from $T$ by removing
${\def\lr#1{\multicolumn{1}{|@{\hspace{.6ex}}c@{\hspace{.6ex}}|}{\raisebox{-.1ex}{$#1$}}}\raisebox{-.6ex}
{$\begin{array}[b]{c}
\cline{1-1}
\lr{\ov{1}}\\
\cline{1-1}
\end{array}$}}$
in the
column  corresponding to the left-most $+$ in $\td{\sigma}$. If there is no such $+$ sign, then we
define $\tf_0 T={\bf 0}$.\vskip 2mm

\noindent \textsc{Case 2.} Suppose that $\epsilon=2$.
Let $T\in \wh{\T}_{n}^{\nw}$ be given.
For $k\geq 1$, let $(t_{2k-1},t_{2k})$ be the pair of entries  in the bottom boxes of the  $(2k-1)$-st and $2k$-th columns of $T$ (from the left), respectively. Note that  $t_{2k-1}$ and $t_{2k}$ are placed in the same row and $t_{2k-1}\geq t_{2k}$. Consider  $\sigma=(\ldots,\sigma_2,\sigma_1)$, where
\begin{equation*}\label{signs}
\sigma_k=
\begin{cases}
- \ , & \text{if  $t_{2k-1},t_{2k}< \ov{1}$ or the $(2k-1)$-st column is empty}, \\
+ \ ,& \text{if $t_{2k-1}=t_{2k}=\ov{1}$}, \\
\ \cdot \ \, ,& \text{otherwise}.
\end{cases}
\end{equation*}
Then we define $\te_0 T$ to be the tableau
obtained from $T$  by adding a domino
${\def\lr#1{\multicolumn{1}{|@{\hspace{.6ex}}c@{\hspace{.6ex}}|}{\raisebox{-.1ex}{$#1$}}}\raisebox{-.6ex}
{$\begin{array}[b]{cc}
\cline{1-1}\cline{2-2}
\lr{\ov{1}}&\lr{\ov{1}}\\
\cline{1-1}\cline{2-2}
\end{array}$}}$
to the bottom of the pair of columns  corresponding to the right-most $-$ in $\td{\sigma}$.
We define $\tf_0 T$ to be the tableau obtained from $T$ by removing a domino
${\def\lr#1{\multicolumn{1}{|@{\hspace{.6ex}}c@{\hspace{.6ex}}|}{\raisebox{-.1ex}{$#1$}}}\raisebox{-.6ex}
{$\begin{array}[b]{cc}
\cline{1-1}\cline{2-2}
\lr{\ov{1}}&\lr{\ov{1}}\\
\cline{1-1}\cline{2-2}
\end{array}$}}$
in the pair of
columns  corresponding to the left-most $+$ in $\td{\sigma}$. If there is no such $+$ sign, then we
define $\tf_0 T={\bf 0}$.

Hence $\wh{\T}_{n}^{\nw}$ is a $U'_q(\wh{\g}_{\wh{I}_n})$-crystal with respect to ${\rm wt}$, $\varepsilon_i$, $\varphi_i$, $\te_i$, $\tf_i$ ($i\in \wh{I}_n$), where
\begin{equation}
\begin{split}
\varphi_0(T)&=\max\left\{\,k\ \Big|\ \tf_0^kT\neq {\bf 0}\,\right\},\ \
\varepsilon_0(T)=\varphi_0(T)- \langle {\rm wt}(T), \wh{h}_0 \rangle.\\
\end{split}
\end{equation}

Define
\begin{equation}
\wh{\kappa}^{\nw} : \wh{\M}_n \longrightarrow \wh{\T}_{n}^{\nw},
\end{equation}
by $\wh{\kappa}^{\nw}(M)={\bf P}(M)^{\nw}$.

By similar arguments as in  \cite[Propositions 3.5 and 6.5]{K11},  we have
\begin{prop}\label{symmetric kappa nw}
$\wh{\kappa}^{\nw}$ is an isomorphism of $U_q({\wh{\g}_{\wh{I}_n}})$-crystals.
\end{prop}

\subsection{KR crystals $\B^{n,s}$}
For $s\geq 1$, let
\begin{equation}
\wh{\M}_n^s=\wh{\M}_n\cap \M_{n\times n}^{\epsilon s}.
\end{equation}
We regard  $\wh{\M}_n^s$ as a subcrystal of $\wh{\M}_n$ and consider a $U'_q(\wh{\g})$-crystal
\begin{equation}
\mathcal{B}^{n,s} =\wh{\M}^s_n\otimes T_{s\wh{\omega}_n}.
\end{equation}

\begin{figure}
\includegraphics[width=6cm, height=14cm]{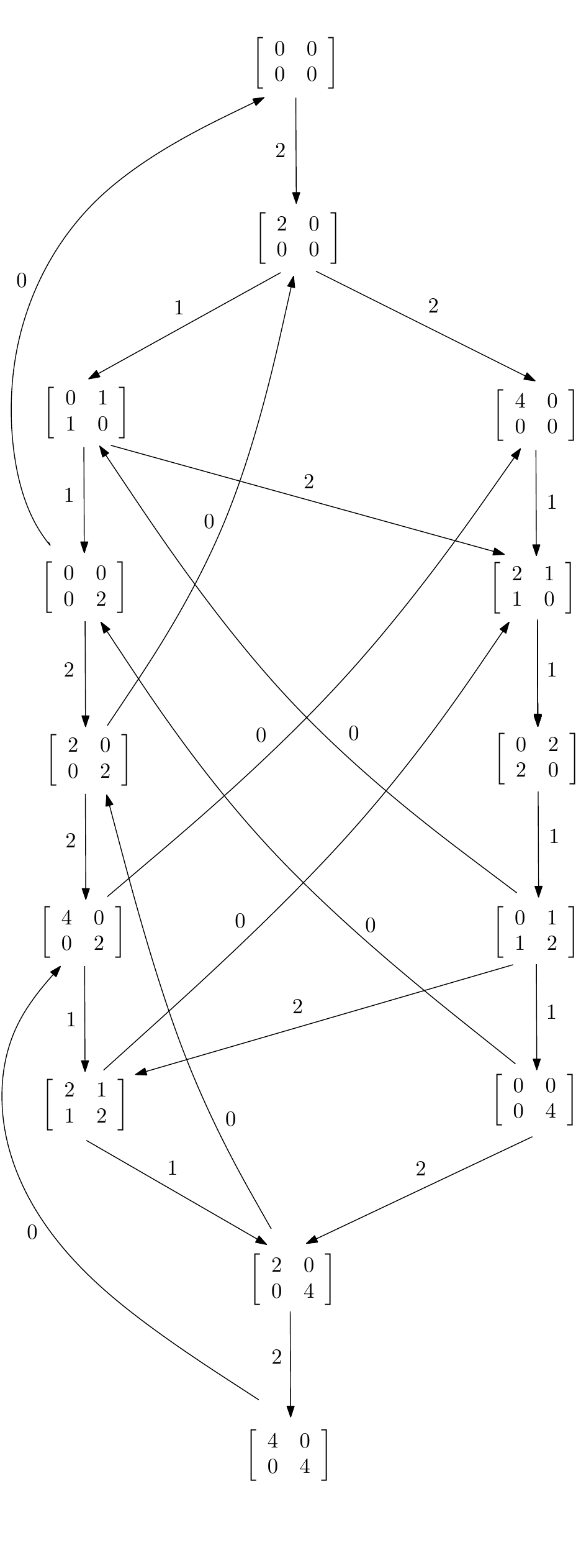}
\caption{The KR crystal graph  $\B^{2,2}$ of type $C_2^{(1)}$.}
\end{figure}

\begin{lem}\label{regularity of folded A type}
$\mathcal{B}^{n,s} $ is a regular $U'_q(\wh{\g})$-crystal that is isomorphic to $\B(s\wh{\omega}_n)$ as a $U_q(\wh{\g}_{\wh{I}_0})$-crystal.
\end{lem}
\pf By (\ref{restricted RSK}), we have bijections
\begin{equation}\label{folded restricted RSK}
\begin{split}
\wh{\kappa}^{\se} : \wh{\M}_n^s &\longrightarrow \wh{\T}_n^{\se,s}, \ \ \ \ \
\wh{\kappa}^{\nw} : \wh{\M}_n^s  \longrightarrow \wh{\T}_n^{\nw,s},
\end{split}
\end{equation}
where
\begin{equation*}
\begin{split}
\wh{\T}_n^{\se,s}&=  \bigsqcup_{\lambda\subset (s^n)}SST_{[\ov{n}]}(\epsilon\lambda^\pi), \ \ \ \ \
\wh{\T}_n^{\nw,s}=  \bigsqcup_{\lambda\subset (s^n)}SST_{[\ov{n}]}(\epsilon\lambda).
\end{split}
\end{equation*}
We may regard $\wh{\T}_n^{\se,s}$ and $\wh{\T}_n^{\nw,s}$ as subcrystals of $\wh{\T}_n^{\se}$ and $\wh{\T}_n^{\nw}$, respectively. Then by Propositions \ref{symmetric kappa se} and \ref{symmetric kappa nw}, the bijections in (\ref{folded restricted RSK}) are isomorphisms of $U_q(\wh{\g}_{\wh{I}_0})$ and $U_q(\wh{\g}_{\wh{I}_n})$-crystals, respectively. On the other hand, by \cite[Remark 5.16]{K09} (or as a special case of \cite[Theorem 6.4] {K11} when $\lambda$ is the empty partition), we have
\begin{equation*}\label{folded covering}
\begin{split}
\mathcal{B}^{n,s}\simeq &\ \wh{\T}_n^{\se,s}\otimes T_{s\wh{\omega}_n} \simeq \B(s\wh{\omega}_n)\, \  \ \ \ \ \ \ \ \ \ \ \ \ \ \ \ \ \ \ \, \text{as a $U_q(\wh{\g}_{\wh{I}_0})$-crystal},\\
\mathcal{B}^{n,s}\simeq &\ \wh{\T}_n^{\nw,s}\otimes T_{s\wh{\omega}_n} \simeq \B(-s\wh{\omega}'_0)\simeq \B(s\wh{\omega}'_0) \, \ \ \ \ \text{as a $U_q(\wh{\g}_{\wh{I}_n})$-crystal}.\\
\end{split}
\end{equation*}
This implies that $\mathcal{B}^{n,s}$ is regular.
\qed\vskip 2mm

Now we have the following, which is the main result in this section.

\begin{thm}\label{main 2}
Let $\B^{n,s}$ be the KR crystal of type $\wh{\g}$ for $s\geq 1$. Then as a $U'_q(\wh{\g})$-crystal, we have
$$\mathcal{B}^{n,s} \simeq \B^{n,s}.$$
\end{thm}
\pf Since $\B^{n,s}\simeq \B(s\wh{\omega}_n)$ as an $U_q(\wh{\g}_{\wh{I}_0})$-crystal (cf.\cite{FOS}), we have $\mathcal{B}^{n,s} \simeq \B^{n,s}$ by Lemmas \ref{ST lemma} and \ref{regularity of folded A type}. \qed

\section{Classically irreducible KR crystals of type $D_{n}^{(1)}$}
\subsection{Affine algebra of type $D_{n}^{(1)}$}
Assume that $\g=D_{n}^{(1)}$ ($n\geq 4$) with $I=\{\,0,1,\ldots,n\,\}$.  Put $I_r=I\setminus\{r\}$ ($r=0,n$), and $I_{0,n}=I_0\cap I_n$. Note that $\g_{I_0}\simeq\g_{I_n}=D_{n}$ and $\g_{I_{0,n}}=A_{n-1}$.

\hskip -3cm
\begin{center}\setlength{\unitlength}{0.16in} \medskip
\begin{picture}(22,4)
\put(1.5,2){\makebox(0,0)[c]{$D_n^{(1)}$ : }}
\put(6,0){\makebox(0,0)[c]{$\bigcirc$}}
\put(6,4){\makebox(0,0)[c]{$\bigcirc$}}
\put(8,2){\makebox(0,0)[c]{$\bigcirc$}}

\put(14,2){\makebox(0,0)[c]{$\bigcirc$}}
\put(16,0){\makebox(0,0)[c]{$\bigcirc$}}
\put(16,4){\makebox(0,0)[c]{$\bigcirc$}}

\put(6.35,0.3){\line(1,1){1.35}} \put(6.35,3.7){\line(1,-1){1.35}}
\put(8.4,2){\line(1,0){1.55}} \put(12,2){\line(1,0){1.55}}

\put(11,1.95){\makebox(0,0)[c]{$\cdots$}}

\put(14.35,2.3){\line(1,1){1.35}} \put(14.35,1.65){\line(1,-1){1.35}}

\put(6,5){\makebox(0,0)[c]{\tiny $\alpha_0$}}
\put(6,-1){\makebox(0,0)[c]{\tiny $\alpha_1$}}
\put(8.2,1){\makebox(0,0)[c]{\tiny $\alpha_2$}}
\put(14,1){\makebox(0,0)[c]{\tiny $\alpha_{n-2}$}}
\put(16,5){\makebox(0,0)[c]{\tiny $\alpha_{n-1}$}}
\put(16,-1){\makebox(0,0)[c]{\tiny $\alpha_{n}$}}
\end{picture}
\end{center}\vskip 5mm

Let
$\epsilon_1=\Lambda_1-\Lambda_0$, $\epsilon_2=\Lambda_2-\Lambda_1-\Lambda_0$, $\epsilon_k=\Lambda_k-\Lambda_{k-1}$ for $k=3,\ldots, n-2$, $\epsilon_{n-1}=\Lambda_{n-1}+\Lambda_n-\Lambda_{n-2}$ and $\epsilon_{n}=\Lambda_{n}-\Lambda_{n-1}$. Then $\bigoplus_{i=1}^n\Z\epsilon_i$ forms a weight lattice of $\g_{I_0}$. Note that $\alpha_i=\epsilon_i-\epsilon_{i+1}$ for $i\in I_{0,n}$, $\alpha_n=\epsilon_{n-1}+\epsilon_n$, and $\alpha_0=-\epsilon_1-\epsilon_2$ in $P_{\rm cl}$.

The fundamental weights for $\g_{I_0}$ are $\omega_i=\sum_{k=1}^i\epsilon_k$ for $i=1,\ldots, n-2$, $\omega_{n-1}=(\epsilon_1+\cdots+\epsilon_{n-1}-\epsilon_n)/2$ and $\omega_{n}=(\epsilon_1+\cdots+\epsilon_{n-1}+\epsilon_n)/2$. We denote the fundamental weights for $\g_{I_n}$ by $\omega'_i$ for $i\in I_n$, where $\omega'_i=\omega_i$ for $i\in I_{0,n}$ and $\omega'_0=-\omega_n$.

\subsection{Young tableaux descriptions of  $\B(s\omega_n)$ and $\B(-s\omega'_0)$ }\label{Crystal of spins}
Consider
\begin{equation}
\begin{split}
{\T}_{n}^{\se}=\bigsqcup_{\substack{\lambda_i': \text{even} \\ \ell(\lambda)\leq n} } SST_{[\ov{n}]}(\lambda^\pi).
\end{split}
\end{equation}
It is a regular $U_q({\g}_{{I}_{0,n}})$-crystal with respect to $\te_i$ and $\tf_i$ for $i\in {I}_{0,n}$, where ${\rm wt}(T)=-\sum_{i\in [n]}m_{\ov{i}}\epsilon_i$ ($m_{\ov{i}}$ is the number of $\ov{i}$'s in $T$) for $T\in {\T}_{n}^{\se}$.

Let us define operators ${\te}_n$, $\tf_n$ on ${\T}_{n}^{\se}$ corresponding to ${\alpha}_n$ as follows:
Let $T\in {\T}_{n}^{\se}$ be given. For $k\geq 1$, let $t_k$  be the entry
in the top of the $k$-th column of $T$ (enumerated from the right). Consider  $\sigma=(\sigma_1,\sigma_2,\ldots)$, where
\begin{equation*}
\sigma_k=
\begin{cases}
+ \ , & \text{if $t_k> \ov{n-1}$ or the $k$-th column is empty}, \\
- \ ,& \text{if  the $k$-th column has both $ \ov{n-1}$ and $ \ov{n}$ as its entries}, \\
\,\cdot\ \ , & \text{otherwise}.
\end{cases}
\end{equation*}
We define $\te_n T$ to be the tableau
obtained from $T$  by removing  the domino
{\tiny ${\def\lr#1{\multicolumn{1}{|@{\hspace{.6ex}}c@{\hspace{.6ex}}|}{\raisebox{-.3ex}{$#1$}}}\raisebox{-.6ex}
{$\begin{array}[b]{c}
\cline{1-1}
\lr{ \ov{n}}\\
\cline{1-1}
\lr{ \ov{n-1}}\\
\cline{1-1}
\end{array}$}}$}
in the column  corresponding to the right-most $-$ in  $\td{\sigma}$.  If there is no such $-$ sign, then we define ${\te}_n T={\bf 0}$.
We define ${\tf}_n T$ to be the tableau  obtained from $T$ by adding a domino
{\tiny ${\def\lr#1{\multicolumn{1}{|@{\hspace{.6ex}}c@{\hspace{.6ex}}|}{\raisebox{-.3ex}{$#1$}}}\raisebox{-.6ex}
{$\begin{array}[b]{c}
\cline{1-1}
\lr{\ov{n}}\\
\cline{1-1}
\lr{ \ov{n-1}}\\
\cline{1-1}
\end{array}$}}$}
on top of the
column  corresponding to the left-most $+$ in $\td{\sigma}$. 
%If there is no such $+$ sign, then we define ${\tf}_n T={\bf 0}$.

Then $\T^{\se}_n$ is a $U_q(\g_{I_0})$-crystal with respect to ${\rm wt}$, $\varepsilon_i, \varphi_i$ and
$\te_i, \tf_i$ $(i\in I_0)$, where 
\begin{equation*}
\varepsilon_n(T)=\max\left\{\,k\ \Big|\ \te_n^k T\neq {\bf 0}\,\right\},\ \
\varphi_n(T)=\varepsilon_n(T)+ \langle {\rm wt}(T), h_n \rangle
\end{equation*}

For $s\geq 1$, let
\begin{equation}\label{se crystal of type Dn1}
\T^{\se,s}_n =\bigsqcup_{\substack{\lambda_i': \text{even} \\ \lambda\subset (s^n)}} SST_{[\ov{n}]}(\lambda^\pi).
\end{equation}
We may consider $\T^{\se,s}_n$ as a subcrystal of $\T^{\se}_n$.
\begin{lem}\label{classical D-1}
 $\T^{\se,s}_n\otimes T_{s\omega_n}$ is isomorphic to $\B(s\omega_n)$ as a $U_q(\g_{I_0})$-crystal.
\end{lem}
\pf First we prove the case when $s=1$. Recall that $\B(\omega_n)$ is  the crystal base of the spin representation of $U_q(\g_{I_0})$, and by \cite{KN} it can be identified with
\begin{equation*}
\{\,v=(i_1,\ldots,i_n)\,|\,i_k=\pm 1,\, i_1\cdots i_n=1 \,\},
\end{equation*}
where ${\rm wt}(v)=\tfrac{1}{2}\sum_{k=1}^n i_k\epsilon_k$ and
\begin{equation*}
\begin{split}
\te_k v=
\begin{cases}
(\ldots,-i_k,-i_{k+1},\ldots), & \text{if $k\in I_{0,n}$  and $(i_k,i_{k+1})=(-1,1)$,} \\
(\ldots,-i_{n-1},-i_{n}),  & \text{if $k=n$ and $(i_{n-1},i_n)=(-1,-1)$},\\
{\bf 0}, & \text{otherwise},
 \end{cases}
\end{split}
\end{equation*}
\begin{equation*}
\begin{split}
\tf_k v=
\begin{cases}
(\ldots,-i_k,-i_{k+1},\ldots), &\text{if $k\in I_{0,n}$  and $(i_k,i_{k+1})=(1,-1)$,} \\
(\ldots,-i_{n-1},-i_{n}),  & \text{if $k=n$ and $(i_{n-1},i_n)=(1,1)$},\\
{\bf 0}, & \text{otherwise}.
 \end{cases}
\end{split}
\end{equation*}
Note that $\T^{\se,1}_n$ is the set of semistandard tableaux with a single column of  even length  no more than $n$.
Define a map $\rho: \T^{\se,1}_n\otimes T_{\omega_n} \rightarrow \B(\omega_n)$ by $\rho(T\otimes t_{\omega_n})=(i_1,\ldots,i_n)$, where $i_k=-1$ if and only if $\ov{k}$ appears in $T$. Note that the empty tableau is mapped to $(1,\ldots,1)$ of weight $\omega_n$. Then it is straightforward to check that  $\rho$ is an isomorphism of $U_q(\g_{I_0})$-crystals.

For $s\geq 1$, consider the map
\begin{equation*}
\iota_s : \T^{\se,s}_n\otimes T_{s\omega_n} \longrightarrow \left(\T^{\se,1}_n \right)^{\otimes s}\otimes T_{s\omega_n} \simeq\left(\T^{\se,1}_n\otimes T_{\omega_n}  \right)^{\otimes s}\simeq \B(\omega_n)^{\otimes s},
\end{equation*}
where for  $\iota_s(T\otimes t_{s\omega_n})=T^1\otimes\cdots \otimes T^s\otimes  t_{s\omega_n}$ ($T^i$ is the $i$-th column of $T$ from the right). Then it is straightforward to check that $\iota_s$ is a strict embedding of $U_q(\g_{I_0})$-crystals, and its image is isomorphic to the connected component of $\emptyset^{\otimes s}\otimes t_{s\omega_n}$, where $\emptyset$ is the empty tableau. Since  $\emptyset^{\otimes s}\otimes t_{s\omega_n}$ is a highest weight element of weight $s\omega_n$ in $\B(\omega_n)^{\otimes s}$, $\T^{\se,s}_n\otimes T_{s\omega_n}$ is isomorphic to $\B(s\omega_n)$.
\qed

Next, consider
\begin{equation}
\ \ \ \ \ {\T}_{n}^{^{\nwarrow}}=\bigsqcup_{\substack{\lambda_i': \text{even} \\ \ell(\lambda)\leq n}} SST_{[\ov{n}]}(\lambda).
\end{equation}
As in \eqref{se crystal of type Dn1}, it is a regular $U_q(\g_{I_{0,n}})$-crystal.
Let us define operators ${\te}_0$, ${\tf}_0$ on ${\T}_{n}^{\nw}$ corresponding to ${\alpha}_0$ as follows:
Let  $T\in {\T}_{n}^{\nw}$ be given.
For $k\geq 1$, let $t_k$  be the entry
in the bottom of the $k$-th column of $T$ (enumerated from the right). Consider  $\sigma=(\ldots,\sigma_2,\sigma_1)$, where
\begin{equation*}
\sigma_k=
\begin{cases}
- \ , & \text{if $t_k< \ov{2}$ or the $k$-th column is empty}, \\
+ \ ,& \text{if  the $k$-th column has both $\ov{1}$ and $\ov{2}$ as its entries,}\\
\, \cdot \ \ , & \text{otherwise}.
\end{cases}
\end{equation*}
We define ${\te}_0 T$ to be the tableau
obtained from $T$  by adding \ 
{\tiny ${\def\lr#1{\multicolumn{1}{|@{\hspace{.6ex}}c@{\hspace{.6ex}}|}{\raisebox{-.3ex}{$#1$}}}\raisebox{-.6ex}
{$\begin{array}[b]{c}
\cline{1-1}
\lr{\ov{2}}\\
\cline{1-1}
\lr{\ov{1}}\\
\cline{1-1}
\end{array}$}}$}\ 
 to the bottom of the column  corresponding to the right-most $-$ in $\td{\sigma}$. 
We define ${\tf}_0 T$ to be the tableau  obtained from $T$ by removing \ {\tiny
${\def\lr#1{\multicolumn{1}{|@{\hspace{.6ex}}c@{\hspace{.6ex}}|}{\raisebox{-.3ex}{$#1$}}}\raisebox{-.6ex}
{$\begin{array}[b]{c}
\cline{1-1}
\lr{\ov{2}}\\
\cline{1-1}
\lr{\ov{1}}\\
\cline{1-1}
\end{array}$}}$} \ 
in the
column  corresponding to the left-most $+$ in $\td{\sigma}$. If there is no such $+$ sign, then we
define ${\tf}_0 T={\bf 0}$. 

Then $\T^{\nw}$ is a $U_q(\g_{I_n})$-crystal with respect to ${\rm wt}$, $\varepsilon_i, \varphi_i$ and
$\te_i, \tf_i$ $(i\in I_n)$, where 
\begin{equation*}
\varphi_0(T)=\max\left\{\,k\ \Big|\ \tf_0^k T\neq {\bf 0}\,\right\},\ \
\varepsilon_0(T)=\varphi_0(T)- \langle {\rm wt}(T), h_0 \rangle
\end{equation*}

For $s\geq 1$, let
\begin{equation}
\T^{\nw,s}_n =\bigsqcup_{\substack{\lambda_i': \text{even} \\ \lambda\subset (s^n)}} SST_{[\ov{n}]}(\lambda).
\end{equation}
We may consider $\T^{\nw,s}_n$ as a subcrystal of $\T^{\nw}_n$.
\begin{lem}\label{classical D-2}
$\T^{\nw,s}_n\otimes T_{s\omega_n}$ is isomorphic to $\B(-s\omega'_0)$ as a $U_q(\g_{I_n})$-crystal.
\end{lem}
\pf The proof is similar to that of Lemma \ref{classical D-1}. Recall that $\B(\omega'_0)$ is also the crystal base of the spin representation of $U_q(\g_{I_0})$, and by \cite{KN} $\B(-\omega'_0)$ can be identified with
\begin{equation*}
\{\,v=(i_1,\ldots,i_n)\,|\,i_k=\pm 1,\, i_1\cdots i_n=1 \,\},
\end{equation*}
where ${\rm wt}(v)=\tfrac{1}{2}\sum_{k=1}^n i_k\epsilon_k$ and
\begin{equation*}
\begin{split}
\te_k v=
\begin{cases}
(\ldots,-i_k,-i_{k+1},\ldots), & \text{if $k\in I_{0,n}$  and $(i_k,i_{k+1})=(-1,1)$,} \\
( -i_{1},-i_{2},\ldots),  & \text{if $k=0$ and $(i_{1},i_2)=(1,1)$},\\
{\bf 0}, & \text{otherwise},
 \end{cases}
\end{split}
\end{equation*}
\begin{equation*}
\begin{split}
\tf_k v=
\begin{cases}
(\ldots,-i_k,-i_{k+1},\ldots), &\text{if $k\in I_{0,n}$  and $(i_k,i_{k+1})=(1,-1)$,} \\
( -i_{1},-i_{2},\ldots),   & \text{if $k=0$ and $(i_{1},i_2)=(-1,-1)$},\\
{\bf 0}, & \text{otherwise}.
 \end{cases}
\end{split}
\end{equation*}
Note that $\T^{\nw,1}_n$ is the set of semistandard tableaux with a single column of even length no more than $n$.
Define a map $\rho: \T^{\nw,1}_n\otimes T_{\omega_n} \rightarrow \B(-\omega'_0)$ by $\rho(T\otimes t_{\omega_n})=(i_1,\ldots,i_n)$, where $i_k=-1$ if and only if $\ov{k}$ appears in $T$. Note that the empty tableau is mapped to $(1,\ldots,1)$ of weight $\omega_n=-\omega'_0$. Then we can check that  $\rho$ is an isomorphism of $U_q(\g_{I_n})$-crystals.

For $s\geq 1$, consider the map
\begin{equation*}
\iota_s : \T^{\nw,s}_n\otimes T_{s\omega_n} \longrightarrow \left(\T^{\nw,1}_n \right)^{\otimes s}\otimes T_{s\omega_n} \simeq\left(\T^{\nw,1}_n\otimes T_{\omega_n}  \right)^{\otimes s}\simeq \B(-\omega'_0)^{\otimes s},
\end{equation*}
where for  $\iota_s(T\otimes t_{s\omega_n})=T^1\otimes\cdots \otimes T^s\otimes  t_{s\omega_n}$ ($T^i$ is the $i$-th column of $T$ from the right). The $\iota_s$ is a strict embedding of $U_q(\g_{I_n})$-crystals, and its image is isomorphic to the connected component of $\emptyset^{\otimes s}\otimes t_{s\omega_n}$, where $\emptyset$ is the empty tableau. Since  $\emptyset^{\otimes s}\otimes t_{s\omega_n}$ is a lowest weight element of weight $s\omega_n=-s\omega'_0$ in $\B(-\omega'_0)^{\otimes s}$, $\T^{\nw,s}_n\otimes T_{s\omega_n}$ is isomorphic to $\B(-s\omega'_0)$.
\qed

\subsection{KR crystals $\B^{n,s}$}\label{KR Bns of type D}

For a semistandard tableau $T$ of skew shape, let $[T]$ denote the equivalence class of $T$ with respect to Knuth equivalence. For $n\geq 4$, let
\begin{equation}
\T_n=\{\,[T]\,|\,T\in \T_n^{\se}\,\}=\{\,[T]\,|\,T\in \T_n^{\nw}\,\}.
\end{equation}
Recall that   under $\te_i$ and $\tf_i$ for $i\in I_{0,n}$, any $T' \in [T]$ generates the same crystal as $T$. Hence, $\mathcal{T}_n$ has a well-defined $U_q(\g_{I_{0,n}})$-crystal structure. Now, for $i=0,n$ and $x=e,f$, we define
\begin{equation}
\begin{split}
\td{x}_i[T]=
\begin{cases}
[\td{x}_0 T^{\nw}], & \text{if $i=0$}, \\
[\td{x}_n T^{\se}], & \text{if $i=n$},
\end{cases}
\end{split}
\end{equation}
where we assume that $[{\bf 0}]={\bf 0}$. Put {\allowdisplaybreaks
\begin{equation}
\begin{split}
{\rm wt}([T])&= {\rm wt}(T), \\
\varepsilon_i([T])&=\varepsilon_i(T), \ \ \ \ \ \ \ \varphi_i([T])=\varphi_i(T) \ \ \ \ (i\in I_{0,n}), \\
\varepsilon_n([T])&=\varepsilon_n(T^{\se}), \ \ \ \  \varphi_n([T])=\varphi_n(T^{\se}),\\
\varepsilon_0([T])&=\varepsilon_n(T^{\nw}), \ \ \ \  \varphi_0([T])=\varphi_n(T^{\nw}).
\end{split}
\end{equation}}
Then, $\T_n$  is a $U'_q(\g)$-crystal with respect to ${\rm wt}$, $\varepsilon_i, \varphi_i$ and
$\te_i, \tf_i$ $(i\in I)$.

Now, for $s\geq 1$, we put $\T_n^s=\{\,[T]\,|\,T\in \T_n^{\se,s}\,\}=\{\,[T]\,|\,T\in \T_n^{\nw,s}\,\}$, which is a subcrystal of $\T_n$, and then define
\begin{equation}
\mathcal{B}^{n,s}=\T_{n}^s\otimes T_{s\omega_n}.
\end{equation}

\begin{lem}\label{regularity of type D}
$\mathcal{B}^{n,s}$  is a regular $U'_q(\g)$-crystal that is isomorphic to $\B(s\omega_n)$ as a  $U_q(\g_{I_0})$-crystal.
\end{lem}
\pf By definition of $\mathcal{B}^{n,s}$ and Lemmas \ref{classical D-1} and \ref{classical D-2}, we have 
\begin{equation*}\label{projections of type Dn1}
\begin{split}
\mathcal{B}^{n,s}&\simeq \T_n^{\se,s}\otimes T_{s\omega_n}\simeq \B(s\omega_n) \,\,\ \text{ as a $U_q(\g_{I_0})$-crystal}, \\
\mathcal{B}^{n,s}&\simeq \T_n^{\nw,s}\otimes T_{s\omega_n}\simeq \B(-s\omega_0') \text{ as a $U_q(\g_{I_n})$-crystal}.
\end{split}
\end{equation*}
This implies that $\mathcal{B}^{n,s}$ is regular.
\qed\vskip 2mm

\begin{thm}\label{main 3}
Let $\B^{n,s}$ be the KR crystal of type ${\g}=D_n^{(1)}$ for $s\geq 1$. Then as a $U'_q({\g})$-crystal, we have
$$\mathcal{B}^{n,s} \simeq \B^{n,s}.$$
\end{thm}
\pf Since $\B^{n,s}\simeq \B(s{\omega}_n)$ as a $U_q(\g_{I_0})$-crystal (cf.\cite{FOS}), 
we have $\mathcal{B}^{n,s} \simeq \B^{n,s}$ by Lemmas \ref{ST lemma} and \ref{regularity of type D}.
\qed

\begin{rem}{\rm
One may expect a matrix realization of $\B^{n,s}$ as in the cases of $A_{n-1}^{(1)}$, $D_{n+1}^{(2)}$ and $C_n^{(1)}$. In fact, there is a variation of RSK map which is a bijection from $\T_n$ to a set of symmetric non-negative integral matrices with trace zero and also an isomorphism of $U_q(A_{n-1})$-crystals (see \cite[Proposition 3.13]{K07} when $m=0$). But there does not seem to be a natural extension to an isomorphism of $U_q(D_n)$-crystals (and hence $U_q(D_n^{(1)})$-crystals). 
}
\end{rem}

\begin{figure}
\includegraphics[width=13cm, height=17cm]{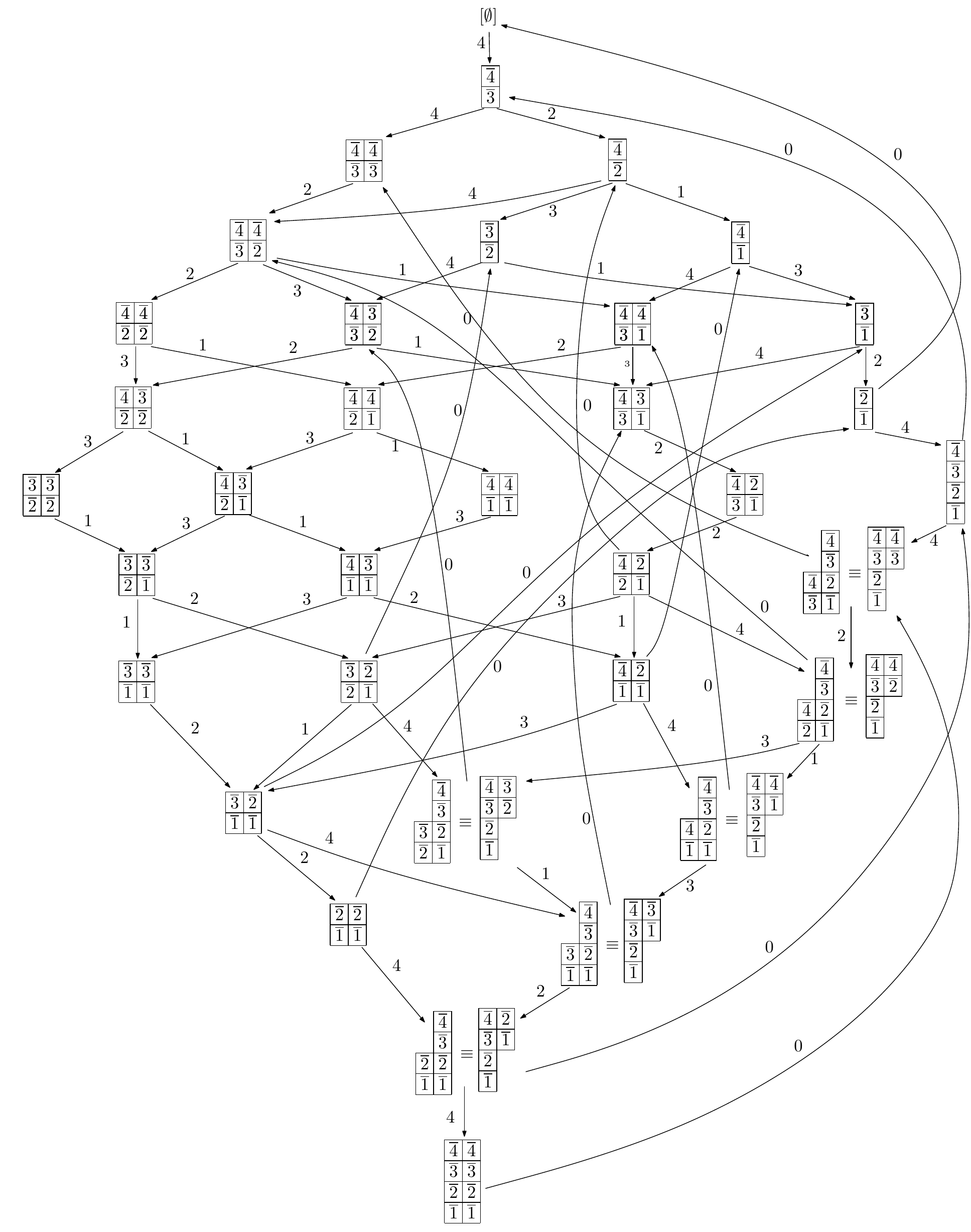}
\caption{The KR crystal graph  $\B^{4,2}$ of type $D_4^{(1)}$. Here $\equiv$ denotes the Knuth equivalence or $U_q(A_3)$-crystal equivalence.}
\end{figure}

\subsection{KR crystals $\B^{n-1,s}$}
Let us give a combinatorial description of $\B^{n-1,s}$ to complete the list of KR crystals associated to exceptional nodes in the Dynkin diagram of classical affine type. Since its construction is very similar to that of $\B^{n,s}$ in Section \ref{KR Bns of type D}, we leave detailed proofs to the reader.

Let
\begin{equation}
\begin{split}
\widetilde{\T}_{n}^{\se}&=\bigsqcup_{\substack{\lambda'_i : \text{odd} \\ \ell(\lambda)\leq n}} SST_{[\ov{n}]}(\lambda^\pi), \ \ \ \ \ 
\widetilde{\T}_{n}^{\nw}=\bigsqcup_{\substack{\lambda'_i : \text{odd} \\ \ell(\lambda)\leq n}} SST_{[\ov{n}]}(\lambda).
\end{split}
\end{equation}
Then $\widetilde{\T}_{n}^{\se}$ (resp. $\widetilde{\T}_{n}^{\nw}$) is a  $U_q(\g_{I_0})$ (resp. $U_q(\g_{I_n})$)-crystal, where $\te_i$, $\tf_i$ for $i\in I_0$ (resp. $i\in I_n$) are defined in the same way as ${\T}_{n}^{\se}$ and ${\T}_{n}^{\nw}$ (Section \ref{Crystal of spins}). For $s\geq 1$, consider subcrystals
\begin{equation}
\begin{split}
\widetilde{\T}_{n}^{\se,s}&=\bigsqcup_{\substack{\lambda'_i : \text{odd} \\ \lambda\subset (s^n)}} SST_{[\ov{n}]}(\lambda^\pi), \ \ \ \ \ 
\widetilde{\T}_{n}^{\nw,s}=\bigsqcup_{\substack{\lambda'_i : \text{odd} \\ \lambda\subset (s^n)}} SST_{[\ov{n}]}(\lambda).
\end{split}
\end{equation}
We have
\begin{equation}
\begin{split}
\widetilde{\T}_{n}^{\se,s}\otimes T_{s\omega_n}&\simeq \B(s\omega_{n-1}) \ \ \ \ \ \ \ \ \text{as a $U_q(\g_{I_0})$-crystal}, \\  
\widetilde{\T}_{n}^{\nw,s}\otimes T_{s\omega_n}&\simeq \B(-s\omega'_1) \ \ \  \ \ \ \ \ \ \text{as a $U_q(\g_{I_n})$-crystal}.
\end{split}
\end{equation}
Let
\begin{equation}
\widetilde{\T}_n=\{\,[T]\,|\,T\in \widetilde{\T}_n^{\se}\,\}=\{\,[T]\,|\,T\in \widetilde{\T}_n^{\nw}\,\},
\end{equation}
where a $U_q'(\g)$-crystal structure on $\widetilde{\T}_n$ is defined in the same way as in ${\T}_n$.

Now, we define
\begin{equation}
\mathcal{B}^{n-1,s}=\widetilde{\T}_{n}^s\otimes T_{s\omega_n}.
\end{equation}
Here $\widetilde{\T}_n^s=\{\,[T]\,|\,T\in \widetilde{\T}_n^{\se,s}\,\}=\{\,[T]\,|\,T\in \widetilde{\T}_n^{\nw,s}\,\}$, which is a subcrystal of $\widetilde{\T}_n$. By the same argument as in Section \ref{KR Bns of type D}, we conclude that
\begin{equation}
\mathcal{B}^{n-1,s}\simeq \B^{n-1,s},
\end{equation}
where $\B^{n-1,s}$ is the KR crystal isomorphic to $\B(s\omega_{n-1})$ as a $U_q(\g_{I_0})$-crystal.

\section{Remarks on $\te_0$ and $\tf_0$}

\subsection{Lusztig involution}
Let $\eta$ be an involutive automorphism of $U_q(A_{n-1})$ given by
\begin{equation*}
\eta(e_i) = f_{n-i}, \ \ \ \eta(f_i)=e_{n-i}, \ \ \ \eta(q^{h_i})=q^{-h_{n-i}} \ \ \ (i=1,\ldots, n-1).
\end{equation*}
Let $w_0$ be the longest element in the Weyl group of $A_{n-1}$. Recall that $w_0(\alpha_i)=-\alpha_{n-i}$ for $i=1,\ldots, n-1$. Let $V$ be a finite dimensional $U_q(A_{n-1})$-module with a crystal base $B$. Then by \cite[Proposition 21.1.2]{Lu}, we have an induced map 
\begin{equation}\label{involution eta}
\eta : B \longrightarrow B
\end{equation}
such that for $b\in B$
\begin{itemize}
\item[(1)] $\eta^2(b) = b$,

\item[(2)] ${\rm wt}(\eta(b)) = w_0({\rm wt}(b))$,

\item[(3)] $\eta(\te_i(b)) =\tf_{n-i}(\eta(b))$ and $\eta(\tf_ib)=\te_{n-i}(\eta(b))$ for $i=1,\ldots, n-1$.
\end{itemize}
In \cite{BZ}, it is shown that  $\eta$ coincides with the Sch\"{u}tzenberger's involution (see e.g. \cite{Fu}) when $B=SST_{[n]}(\lambda)$ for $\lambda\in \cP$ with $\ell(\lambda)\leq n$. Indeed, for $T\in SST_{[n]}(\lambda)$, let $T'$ be the tableau obtained  by $180^{\circ}$-rotation of $T$ and replacing $i$ with $n-i+1$. Then $\eta(T)$ is the unique tableau in $SST_{[n]}(\lambda)$ such that $\eta(T)$ is Knuth equivalent to $T'$.

Similarly, one can define $\eta$ on a crystal associated to a finite dimensional representation of $U_q(A_{m-1}\oplus A_{n-1})$ for $m,n\geq 2$. 
          
Based on our combinatorial descriptions, we have  the following characterization of  $\te_0$ and $\tf_0$ on classically irreducible KR crystals in terms of an involution \eqref{involution eta} on an underlying classical crystal of type $A$.
 
\begin{prop}\label{characterization of x_0}
Let $\B^{r,s}$  be a classically irreducible KR crystal of type $\g$ $(s\geq 1)$ $($that is, for $r=1,\ldots, n-1$ when $\g=A_{n-1}^{(1)}$, $r=n$ when $\g=D_{n+1}^{(2)}$, $C_n^{(1)}$,  $r=n,n-1$ when $\g=D_n^{(1)}$, and  $s\geq 1$$)$. Let $\eta$ denote the involution \eqref{involution eta} on $\B^{r,s}$ as a crystal of type $\g_{J}$ with $J=I\setminus \{0,r\}$. Then we have on $\B^{r,s}$
\begin{equation*}
\begin{split}
\te_0=\eta\circ \tf_r \circ \eta, \ \ \ \
\tf_0=\eta\circ \te_r \circ \eta.
\end{split}
\end{equation*} 
\end{prop}
\pf Throughout the proof,  we assume that $x=e$ (resp. $f$) when $y=f$ (resp. $e$).

\noindent\textsc{Case 1}. $\B^{r,s}$ of type $A_{n-1}^{(1)}$ for $r=1,\ldots, n-1$ and $s\geq 1$.
Note that $\g_J\simeq A_{r-1}\oplus A_{n-r-1}$ in this case.
Consider  
\begin{equation*}
\pi : \M_{r\times (n-r)} \longrightarrow \M_{r\times (n-r)},
\end{equation*}
where $\pi(M)$ is the matrix obtained by $180^\circ$ rotation of $M$. Note that $\pi^2={\rm id}$.
By definition of $\te_0$ and $\tf_0$ on $\M_{r\times (n-r)}$, we have
\begin{equation}\label{symmetry of e_0 and f_0}
\widetilde{x}_0=\pi\circ \widetilde{y}_r \circ \pi.
\end{equation}

Now, let $M=M({\ba},{\bb})$ be given with $\ba=\ov{i_1}\ldots \ov{i_k}$. Then $\pi(M)=M(\ba^\pi,\bb^\pi)$ with $\ba^\pi=\ov{r-i_k+1}\ldots \ov{r-i_1+1}$. Also, if $M^t=M({\bc},{\bd})$ with $\bc=j_1\ldots j_l$, then $\pi(M^t)=M(\bc^\pi,\bd^\pi)$ with $\bc^\pi=(n-j_l+r+1)\ldots (n-j_1+r+1)$. This implies that 
\begin{equation}\label{involution pi}
\begin{split}
\widetilde{x}_i M \neq {\bf 0}\ \ \ \ &\Longleftrightarrow\ \ \ \ \widetilde{y}_{r-i}(\pi(M))\neq {\bf 0} \ \ \  \ \ \ (i=1,\ldots,r-1),\\
\widetilde{x}_i M \neq {\bf 0}\ \ \ \ &\Longleftrightarrow\ \ \ \ \widetilde{y}_{n-i+r}(\pi(M))\neq {\bf 0} \ \ \ (i=r+1,\ldots,n-1).
\end{split}
\end{equation}

On the other hand, we also have 
\begin{equation}\label{involution S}
\begin{split}
\widetilde{x}_i M \neq {\bf 0}\ \ \ \ &\Longleftrightarrow\ \ \ \ \widetilde{y}_{r-i}(\eta(M))\neq {\bf 0} \ \ \  \ \ \ (i=1,\ldots,r-1),\\
\widetilde{x}_i M \neq {\bf 0}\ \ \ \ &\Longleftrightarrow\ \ \ \ \widetilde{y}_{n-i+r}(\eta(M))\neq {\bf 0} \ \ \ (i=r+1,\ldots,n-1).
\end{split}
\end{equation}

Let $M=(m_{\ov{i} j})$ be a $\g_J$-highest weight element in $\M_{r\times (n-r)}$, where  $m_{\ov{i} j}=0$ unless $i=j$, and $m_{\ov{r}\, r+1}\geq m_{\ov{r-1}\, r+2}\geq m_{\ov{r-2}\, r+2}\geq \ldots$.
Then it is straightforward to check that $\pi(M)=\eta(M)$.
From \eqref{involution pi} and \eqref{involution S}, it follows that $\pi=\eta$ on $\M_{r\times (n-r)}$. Therefore, by \eqref{symmetry of e_0 and f_0}, we have 
\begin{equation}\label{x_0 on type A}
\widetilde{x}_0=\eta\circ \widetilde{y}_r \circ \eta.
\end{equation}
Since $\B^{r,s}$ is a subcrystal of $\M_{r\times (n-r)}\otimes T_{s\omega_r}$, the relation \eqref{x_0 on type A} also holds on $\B^{r,s}$.

\vskip 5mm

\noindent\textsc{Case 2}. $\B^{n,s}$ of type $D_{n+1}^{(2)}$, $C_n^{(1)}$ for $s\geq 1$.

The proof is similar to that of Case 1. Note that $\pi$ on $\M_{2n\times 2n}$ induces an involution on $\wh{\M}_n$, and by \eqref{symmetry of e_0 and f_0}
\begin{equation}
\wh{x}_0=\pi\circ \wh{y}_n \circ \pi.
\end{equation}

Let $M\in \wh{\M}_n$ be given.
By \eqref{involution pi}, we have
\begin{equation}\label{involution pi-2}
\wh{x}_i M \neq {\bf 0}\ \ \ \ \Longleftrightarrow\ \ \ \ \wh{y}_{n-i}(\pi(M))\neq {\bf 0} \ \ \  \ \ \ (i=1,\ldots, n-1).
\end{equation}
On the other hand, since $\g_J\simeq A_{n-1}$, we also have 
\begin{equation}\label{involution S-2}
\wh{x}_i M \neq {\bf 0}\ \ \ \ \Longleftrightarrow\ \ \ \ \wh{y}_{n-i}(\eta(M))\neq {\bf 0} \ \ \  \ \ \ (i=1,\ldots,n-1).
\end{equation}
Since $\pi(M)=\eta(M)$ for a $\g_J$-highest weight element $M$ in $\wh{\M}_n$, we have $\pi=\eta$ on $\wh{\M}_n$ by \eqref{involution pi-2} and \eqref{involution S-2}, which implies that
\begin{equation}\label{x_0 on folded case}
\wh{x}_0=\eta\circ \wh{y}_n \circ \eta.
\end{equation}
Since $\mathcal{B}^{n,s} $ is a subcrystal of $\wh{\M}_n\otimes T_{s\wh{\omega}_n}$, the relation \eqref{x_0 on folded case} also holds on $\B^{n,s}$.
\vskip 5mm

\noindent\textsc{Case 3}. $\B^{r,s}$ of type $D_{n}^{(1)}$ for $r=n,n-1$ and $s\geq 1$.

Let us prove the case $\B^{n,s}$. The proof for $\B^{n-1,s}$ is almost the same and we leave it to the readers.

Let $[T]\in {\T}_n$ be given. Define a map 
\begin{equation}
\pi : {\T}_n \longrightarrow {\T}_n,
\end{equation}
where $\pi([T])=[T']$ and $T'$ is obtained from by $180^\circ$ of $T$ and replacing each entry $\ov{i}$ in $T$ with $\ov{n-i+1}$.
By definition, it is not difficult to see that $\pi^2={\rm id}$ and 
\begin{equation}
\widetilde{x}_i T \neq {\bf 0}\ \ \ \ \Longleftrightarrow\ \ \ \ \widetilde{y}_{n-i}T'\neq {\bf 0} \ \ \  \ \ \ (i=1,\ldots,n-1).
\end{equation}
This implies that $T'=\eta(T)$.

Moreover, if $T$ is of normal shape, then we have by definition of $\widetilde{x}_0$ and $\widetilde{y}_n$ (see Section \ref{Crystal of spins})
\begin{equation}
\widetilde{x}_0([T])=\left(\pi \circ \widetilde{y}_n \circ\pi\right) ([T]).
\end{equation} 
Since the action of $\eta$ is also well-defined on ${\T}_n$ (that is, $\eta([T])=[\eta(T)]$), we conclude that
\begin{equation}\label{x_0 on D_n^1}
\widetilde{x}_0=\eta \circ \widetilde{y}_n \circ \eta.
\end{equation}
Since $\B^{n,s}$ is a subcrystal of $\T_n\otimes T_{s\omega_n}$, the relation \eqref{x_0 on D_n^1} also holds on $\B^{n,s}$.
\qed

\subsection{A connection with the Sch\"{u}tzenberger's promotion operator}
Let ${\rm\bf pr}$ be the Sch\"{u}tzenberger's promotion operator on ${SST}_{[n]}(\lambda)$ ($\lambda\in \cP$) \cite{Sc}, which satisfies for $T\in SST_{[n]}(\lambda)$ with ${\rm wt}(T) = m_1\epsilon_1+m_2\epsilon_2+\cdots + m_n\epsilon_n$
\begin{itemize}
\item[(1)] ${\rm wt}({\rm\bf pr}(T)) = m_n\epsilon_1+m_1\epsilon_2+\cdots + m_{n-1}\epsilon_n$,

\item[(2)] ${\rm\bf pr}(\te_i T) = \te_{i+1}( {\rm\bf pr}(T))$ and ${\rm\bf pr}(\tf_i T) = \tf_{i+1}( {\rm\bf pr}(T))$ for $i=1,\ldots, n-2$.
\end{itemize}
Note that  ${\rm \bf pr}$  is the unique map on $SST_{[n]}(\lambda)$ satisfying  (1) and (2), and  ${\rm \bf pr}$  is of order $n$ if and only if $\lambda$ is a rectangle (see \cite[Proposition 3.2]{ST}).

It is shown  in \cite{S} that on $\B^{r,s}$ of type $A_{n-1}^{(1)}$ ($r=1,\ldots, n-1$, $s\geq 1$)
\begin{equation*}
\te_0 = {\rm\bf pr}^{-1} \circ \te_1 \circ {\rm \bf pr}, \ \ \ \tf_0 = {\rm\bf pr}^{-1} \circ \tf_1 \circ {\rm \bf pr}.
\end{equation*}

Suppose that $\g=A_{n-1}^{(1)}$. For $k\in I$, let $\eta_k$ denote the involution \eqref{involution eta} on crystals of type $\g_{I_{0,k}}$. Here $\g_{I_{0,0}}=\g_{I_0}$.
Let $\lambda\in \cP$ be given with $\ell(\lambda)\leq n$. 
 Put $\xi=\eta_1\circ\eta_0$. By definition of $\xi$, it is straightforward to check that 
\begin{itemize}
\item[(1)] ${\rm wt}(\xi(T)) = m_n\epsilon_1+m_1\epsilon_2+\cdots + m_{n-1}\epsilon_n$,

\item[(2)] $\xi(\te_i T) = \te_{i+1}( \xi(T))$ and $\xi(\tf_i T) = \tf_{i+1}(\xi(T))$ for $i=1,\ldots, n-2$.
\end{itemize}
By the uniqueness of ${\rm \bf pr}$, we have  on $SST_{[n]}(\lambda)$
\begin{equation}
{\rm\bf pr} = \eta_1\circ \eta_0.
\end{equation} 

\begin{lem}\label{eta_0 e eta_0} For $b\in \B^{r,s}$, we have  $\eta_0(\te_0 b) =\tf_0 (\eta_0(b))$ and $\eta_0(\tf_0 b) =\te_0 (\eta_0(b))$.
\end{lem}
\pf First, we claim that 
\begin{equation}\label{eta e eta}
\te_0 =\eta_1\circ\tf_1\circ\eta_1, \ \ \ \ \tf_0 =\eta_1\circ\te_1\circ\eta_1.
\end{equation}
Note that ${\rm \bf pr}^n={\rm id}_{\B^{r,s}}$. We have
\begin{equation*}
\begin{split}
{\rm \bf pr}\circ\te_{n-1} &= {\rm \bf pr}^{n-1}\circ \te_{1} \circ {\rm \bf pr}^{-n+2} ={\rm \bf pr}^{-1}\circ \te_{1} \circ {\rm \bf pr}^{2} = \te_0 \circ {\rm \bf pr}.
\end{split}
\end{equation*}
Since ${\rm \bf pr}=\eta_1\circ\eta_0$, we have
\begin{equation*}
\begin{split}
\te_0=\eta_1\circ\eta_0\circ\te_{n-1}\circ\eta_0\circ\eta_1 &= \eta_1\circ\eta_0\circ\eta_0\circ \tf_1\circ\eta_1= \eta_1\circ\tf_1\circ\eta_1.
\end{split}
\end{equation*}
Similarly, we have $\tf_0 =\eta_1\circ\te_1\circ\eta_1$.
Now, by \eqref{eta e eta}, we have
\begin{equation*}
\begin{split}
\eta_0\circ\te_{0} &= \eta_0\circ{\rm \bf pr}^{-1}\circ \te_{1} \circ {\rm \bf pr} =\eta_0\circ\eta_0\circ\eta_1\circ \te_{1} \circ \eta_1\circ\eta_0=\tf_0 \circ \eta_0.
\end{split}
\end{equation*}
Since $\eta$ is an involution, we also have $\eta_0\circ \tf_0= \te_0\circ\eta_0$.
\qed

\begin{prop}\label{pr^k} Let $\B^{r,s}$ be a KR crystal of type $A_{n-1}^{(1)}$ for $1\leq r\leq n-1$ and $s\geq 1$. Then we have 
$$
{\rm \bf pr}^k = \eta_k \circ \eta_0,
$$
on $\B^{r,s}$ for $1\leq k \leq n-1$.
\end{prop}
\pf It is not difficult to see that the highest (resp. lowest) weight elements in $\B^{r,s}$ as a $U_q(\g_{I_{0,k}})$-crystal are parametrized by the partitions  $\lambda\subset (s^r)$, say $b^{\rm h.w.}_\lambda$ (resp. $b^{\rm l.w.}_\lambda$). 
Note that $\eta_k \circ \widetilde{x}_i = \widetilde{y}_{n+k-i} \circ \eta_k$ for $i\in I_{0,k}$ and $\eta_k (b_\lambda^{\rm h.w.})=b_\lambda^{\rm l.w.}$ for $\lambda\subset (s^r)$. Here we assume that $x=e$ (resp. $f$) when $y=f$ (resp. $e$), and the indices are assumed to be the elements of $\mathbb{Z}_n$.

Let $\xi_k = {\rm \bf pr}^k\circ \eta_0$. It suffices to show that $\xi_k=\eta_k$. First, it is straightforward to check that
\begin{equation}\label{xi_k}
\begin{split}
\xi_k \circ \widetilde{x}_i &= \widetilde{y}_{n+k-i} \circ \xi_k, \\
\end{split}
\end{equation}
for $i\in I_{0,k}$. Since  $\xi_k (b_\lambda^{\rm h.w.})$ is a lowest weight element as a $U_q(\g_{I_{0,k}})$-crystal by \eqref{xi_k}, and ${\rm wt}(\xi_k (b_\lambda^{\rm h.w.}))={\rm wt}(b_\lambda^{\rm l.w.})$, we have $\xi_k (b_\lambda^{\rm h.w.})=b_\lambda^{\rm l.w.}$.  Therefore $\xi_k(b) =\eta_k(b)$ for  $b\in \B^{r,s}$.
\qed

\begin{cor}\label{eta k} Under the above hypothesis, we have 
\begin{equation*}
\begin{split}
\te_0=\eta_k\circ \tf_k \circ \eta_k, \ \ \ \
\tf_0=\eta_k\circ \te_k \circ \eta_k,
\end{split}
\end{equation*} 
on $\B^{r,s}$ for $1\leq k\leq n-1$.
\end{cor}
\pf Since ${\rm \bf pr}^{-k}\circ \te_k \circ {\rm \bf pr}^{k}=\te_0$, we have by Proposition \ref{pr^k}
\begin{equation*}
\begin{split}
\eta_0\circ \eta_k\circ \te_k\circ \eta_k\circ \eta_0 = \te_0
\end{split}
\end{equation*} 
Hence, we have $\eta_k\circ \te_k\circ \eta_k=\eta_0\circ\te_0\circ\eta_0=\tf_0$ by Lemma \ref{eta_0 e eta_0}. Similarly, we have $\tf_0=\eta_k\circ \te_k \circ \eta_k$.
\qed

\begin{rem}{\rm By Proposition \ref{pr^k}   the operators of $\eta_1$ and $\eta_0$ on $\B^{r,s}$ generate the Dihedral group of order $2n$.  When $k=r$, Corollary \ref{eta k} also implies Proposition \ref{characterization of x_0} for type $A_{n-1}^{(1)}$.

}
\end{rem}


{\small
\begin{thebibliography}{K09}
\bibitem{BZ}
A. Berenstein, A. Zelevinsky, {\em Canonical bases for the quantum group of type $A_r$ and piecewise-linear combinatorics}, Duke Math. J. \textbf{82} (1996) 473--502.

\bibitem{DK}
V.~I. Danilov, G.~A. Koshevoy, {\em Bi-crystals and crystal
$(GL(V),GL(W))$ duality}, RIMS preprint, (2004) no. 1458.

\bibitem{Fu}
W. Fulton, {\em Young tableaux}, London Mathematical Society Student
Texts, 35. Cambridge University Press, Cambridge, 1997.

\bibitem{FOS}
G. Fourier, M. Okado, A. Schilling, {\em Kirillov-Reshetikhin crystals for nonexceptional types}, Adv. Math. \textbf{222} (2009) 1080--1116.

\bibitem{HKOTY}
G. Hatayama, A. Kuniba, M. Okado, T. Takagi, Y. Yamada, {\em Remarks on fermionic formula}, in: Contemp. Math., \textbf{248}, (1999)  243--291.

\bibitem{HKOTT}
G. Hatayama, A. Kuniba, M. Okado, T. Takagi, Z. Tsuboi, {\em Paths, crystals and fermionic formulae}, in: MathPhys Odyssey 2001, in: Prog. Math. Phys., vol. 23, BirkhŠuser Boston, Boston, MA, 2002, 205--272.

\bibitem{HK}
J. Hong, S.-J. Kang, {\em Introduction to Quantum Groups and Crystal Bases},
Graduate Studies in Mathematics Vol. 42, Amer. Math. Soc., Providence, RI,
2002.

\bibitem{K} V.~Kac, {\em Infinite-dimensional Lie algebras}, Third
edition, Cambridge University Press, Cambridge, 1990.

\bibitem{KMN2}
S.-J. Kang, M. Kashiwara, K.C. Misra, T. Miwa, T. Nakashima, A. Nakayashiki, {\em Perfect crystals of quantum affine
Lie algebras}, Duke Math. J. \textbf{68} (1992) 499--607.

\bibitem{Kas91}
M. Kashiwara, \emph{On crystal bases of the $q$-analogue of
universal enveloping algebras}, Duke Math. J. \textbf{63} (1991),
465--516.

\bibitem{Kas94}
M. Kashiwara, \emph{On crystal bases}, Representations of groups,
CMS Conf. Proc., 16, Amer. Math. Soc., Providence, RI, (1995)
155--197.

\bibitem{Kas96}
M. Kashiwara, {\em Similarity of crystal bases}, Contemp. Math.
\textbf{194} (1996), 177--186.

\bibitem{KN}
M. Kashiwara, T. Nakashima, {\em Crystal graphs for representations
of the $q$-analogue of classical Lie algebras}, J. Algebra
\textbf{165} (1994) 295--345.

\bibitem{KR}
A.N. Kirillov, N.Yu. Reshetikhin, {\em Representations of Yangians
and multiplicities of the inclusions of the irreducible components
of the tensor product of representations of simple Lie algebras.}
Zap. Nauchn. Sem. Leningrad. Otdel. Mat. Inst. Steklov. (LOMI) 160
(1987), Anal. Teor. Chisel i Teor. Funktsii. 8, 211.221, 301;
translation in J. Soviet Math. 52 (1990), no. 3, 3156.3164

\bibitem{K07}
J.-H. Kwon, {\em Crystal graphs for Lie superalgebras and Cauchy
decomposition}, J. Algebraic Combin. \textbf{25} (2007) 57--100.

\bibitem{K09}
J.-H. Kwon, {\em Demazure crystals of generalized Verma modules and
a flagged RSK correspondence}, J. Algebra \textbf{322} (2009), 2150--2179.

\bibitem{K10}
J.-H. Kwon, {\em Crystal bases of modified quantized enveloping algebras and a
double RSK correspondence}, J. Combin. Theory Ser. A. \textbf{118} (2011) 2131--2156.


\bibitem{K11}
J.-H. Kwon {\em Littlewood identity and crystal bases},  preprint, arXiv:1106.5286, to appear in Adv. Math.

\bibitem{La}
A. Lascoux, {\em Double crystal graphs}, Studies in Memory of Issai
Schur, Progress in Math. \textbf{210}, Birkh\"{a}user (2003),
95--114.

\bibitem{Lu}
G.Lusztig, {\em Introduction to quantum groups}, Birkha\"{u}ser, Boston 1993.

\bibitem{Mac95}
I. G. Macdonald, {\em Symmetric functions and Hall polynomials},
Oxford University Press, 2nd ed.,  1995.

\bibitem{OS}
M. Okado, A. Schilling, {\em Existence of Kirillov-Reshetikhin crystals for nonexceptional types}, Represent. Theory \textbf{12}
(2008) 186--207.

\bibitem{Sc}
M.-P. Sch\"{u}tzenberger, {\em Promotion des morphisms dÕensembles ordonnes}, Disc. Math. \textbf{2} (1972) 73--94.

\bibitem{ST}
A. Schilling, P. Tingley, {\em Demazure crystals,
Kirillov-Reshetikhin crystals, and the energy functions},
arXiv:1104.2359.

\bibitem{S}
M. Shimozono, {\em Affine type A crystal structure on tensor products of rectangles, Demazure characters, and nilpotent
varieties}, J. Algebraic Combin. \textbf{15} (2002) 151--187.

\end{thebibliography}}
\end{document}